\documentclass[A4]{amsart}
\usepackage{geometry}                
\usepackage{graphicx}
\usepackage{amssymb}
\usepackage{epstopdf}
\newcommand{\CC}{\hbox{{$\mathcal C$}}}
\newcommand{\CA}{\hbox{{$\mathcal A$}}}

\newcommand{\CL}{\hbox{{$\mathcal L$}}}
\newcommand{\CB}{\hbox{{$\mathcal B$}}}

\newcommand{\CV}{\hbox{{$\mathcal V$}}}

\newcommand{\cg}{\mathfrak{g}}

\newcommand{\C}{\mathbb{C}}
\newcommand{\R}{\mathbb{R}}

\newcommand{\Z}{\mathbb{Z}}

\newcommand{\del}{\partial}

\newcommand{\extd}{{\rm d}}

\newcommand{\isom}{{\cong}}
\newcommand{\eps}{{\epsilon}}
\newcommand{\tens}{\mathop{\otimes}}
\newcommand{\la}{{\triangleright}}

\newcommand{\id}{{\rm id}}

\def\lcross{{>\!\!\!\triangleleft}}

\newcommand{\nquad}{\kern-60pt}

\newcommand{\und}[1]{{\underline {#1}}}
\newcommand{\eqn}[2]{\begin{equation}#2\label{#1}\end{equation}}

\begin{document}
 
\title[Gauge theory on nonassociative spaces]{\rm \large Gauge theory on nonassociative spaces}
\author{S. Majid}%

\address{School of Mathematical Sciences\\
Queen Mary, University of London\\ 327 Mile End Rd,  London E1
4NS, UK}

\date{14 June 2005}                                           

\begin{abstract}We show how to do gauge theory on  
the octonions and other nonassociative algebras such as `fuzzy $\R^4$' models proposed
in string theory. We use the theory of quasialgebras obtained by cochain twist introduced previously. The gauge theory in this case is twisting-equivalent to usual gauge theory on the underlying classical space. We give a general  $U(1)$-Yang-Mills example for any quasi-algebra and a full description of the moduli space of flat connections in this theory for the cube $\Z_2^3$ and hence for the octonions. We also obtain further results about the octonions themselves; an explicit Moyal-product description of them as a nonassociative quantisation of functions on the cube, and a characterisation of their cochain twist as invariant under Fourier transform.\end{abstract}

\maketitle

\section{Introduction}

There has been a lot of interest recently in `nonassociative geometry' as a further extension
of the ideas of noncommutative geometry, with now the `coordinate algebra' allowed to be
nonassociative. The framework which we use of 'quasialgebras' was already established and used to describe the octonions as 'quasispaces' some years ago \cite{AlbMa:qua}. These were, moreover, constructed as a `cochain twist' of a classical associative space. Differential geometry on such quasispaces was introduced in \cite{AkrMa:bra} and in this paper we add 'gauge theory'. 

The need for nonassociative geometry for noncommutative differential forms (even when the coordinate algebra itself remains associative) was shown in \cite{BegMa:sem}, where it was proven that all differential
form algebras on the standard q-deformation quantum groups, if they are to be bicovariant and to
have classical dimensions, must indeed be nonassociative. Thus the usual assumption in noncommutative geometry, including in \cite{Con}, that differential forms should be associative, appears to be too strong. From a physics point of view also, there are suggestions that the world volume algebras on certain string theories are naturally nonassociative,
and this has been realised quite concretely in some form in the context of reduced matrix models, see \cite{Ram:sph, HoRam:geo, Ram:gau}. In the latter is posed the problem of gauge theory on such spaces, with apparently higher order differentials being required. We start by making precise what is fairly clear that the simplified `fuzzy' algebras in \cite{Ram:gau} are indeed quasialgebras in our required sense. We then show that in this case
there is a natural formulation of gauge theory on them looking much more like the classical case. We describe this theory for any quasialgebra (or algebra in a nonassociative monoidal category) at an algebraic level and give a general construction for examples equivalent to $U(1)$-Yang-Mills in the associative case. The framework allows for nonAbelian gauge theory as well. Also, we do not discuss Lagrangians here but all of the necessary data and methods for these are known in the associative case, see notably \cite{MaRai:ele,Ma:phy}, and apply equivalently to quasialgebras obtained by  cochain twists.

As well as covering the string-motivated example, we explore fully the octonions as finite  'quasi-geometries' par excellence.  We show that the cochain $F(\vec a,\vec b)$ in \cite{AlbMa:qua} that modifies the group algebra of the cube $\Z_2^3$ to the octonion product has the very remarkable feature of being invariant under $\Z_2^3$-Fourier transform. Using this, we also find an explicit  more geometrical  $\bullet$-product  description of the octonion  as a nonassociative quantisation of the coordinate algebra on the Fourier-dual cube $\Z_2^3$ by means if a (finite difference) bidifferential operator. This is in the spirit of the Moyal-product of functions on $\R^n$, but now nonassociative. The associative quantisation (Clifford algebra) case is also covered.   

The paper begins in Section 2 with a brief introduction the theory of quasialgebras obtained by cochain twist\cite{AlbMa:qua, AkrMa:bra}, as algebras in a (symmetric) monoidal category. Sections~2.1 and~2.2 respectively outline the continuous case deforming $\R^n$ and the finite case deforming  group algebras. In Secton 3.1 we recall from \cite{Ma:rem} the formulation of gauge theory in such a general monoidal category and the diagrammatic notation for it. Section~3.2 appiles this at an algebraic level to describe gauge theory on cochain twist quasi-algebras in general. Section~3.3 gives a canonical general example where the 'gauge group' can be chosen canonically. Although appearing nonAbelian (and nonassociative) we show that this particular choice gives a theory equivalent to the undeformed $U(1)$-Yang-Mills theory, Note that in noncommutative geometry even the $U(1)$ theory has $F(\alpha)=\extd\alpha+\alpha\wedge\alpha$ and we use the phrase Yang-MIlls to distinguish this nonlinear theory from the Maxwell case where $F(\alpha)=\extd\alpha$. 

In Section 4 we apply the theory the quasialgebra versions of $\R^n$  of interest in \cite{Ram:gau} under heading of a simplified `fuzzy $\R^n$'. Section~4.1 introduces the required nonassociative differential calculus and Section~4.2 the promised gauge theory. Finally, in Section~5 we apply the theory the octonions. Section~5.1 warms up with the new results about the octonions as $\bullet$-product. Section~5.2 has the gauge theory worked out for the octonions. In fact the example of deformed nonassociative gauge theory that we finally arrive at here takes the remarkably workable form
\[ F_\bullet(\alpha)=\extd\alpha+\sum F(|\alpha|,|\alpha'|)\alpha\bullet\alpha'\]
\[ \alpha_\bullet^\gamma= \sum F(|\gamma^{-1}|,|\gamma|)F(|\alpha|,|\gamma|) (\gamma^{-1}\bullet\alpha)\bullet\gamma+\sum F(|\gamma^{-1}|,|\gamma|)\gamma^{-1}\bullet \extd\gamma\]
where the sum is over the different graded Fourier components of each object and to this end $\alpha'$ denotes a second independent copy of $\alpha$.  Such a description also works  for the fuzzy-$\R^4$ if one works in terms of plane waves and their differentials; this is already the case for the octonions where the generators have in our picture the interpretation of deformed plane waves on the cube. In the octonions case $F(\vec a,\vec b)$ has values $\pm 1$ but is not simply an exponential bilinear of the vector degrees (the 3-momentum) as its exponent has cubic terms. It is not known if such quasi-geometry of the octonions has a direct physical role, but see for example \cite{BHTZ}. We also note the link between the octonions and particle physics\cite{Dix}; their geometry might play a role in the context of the direct product of spacetime by the finite geometry.  

Section~5.3 fills a gap in the literature, namely a complete description of the moduli space of flat $U(1)$-Yang-Mills fields  up to gauge transformation on $\Z_2^2$ and $\Z_2^3$, using the same methods as for the symmetric group $S_3$ in \cite{MaRai:ele}. The above equivalence means that the $\Z_2^3$ case also classifies flat connections in the nonassociative theory on the octonions. We note that flat connections on finite groups are also of interest in pure mathematics in connection with Schubert calculus on flag varieties \cite{Ma:yan}. Going back to physics, the quantum $U(1)$-Yang-Mills theory on $\Z_2^2$ is fully worked out in \cite{Ma:phy} and is renormalisable and computable. The $\Z_2^3$ and octonion cases could in principle be similarly computed. Thus would be one of several directions for further work. 

We also note the related paper \cite{BegMa:qua} where cochain twists are used to describe associative quantisations in which the differential calculus, however, is nonassociative. It turns out that several popular associative quantisations in physics fall into this category; the algebra of coordinates is associative but the nonassociative gauge theory described here still plays a role in view of the differential calculus. Examples in this category include $U(\cg)$ as quantisation of the Kirillov-Kostant bracket, now expressed as a cochain twist at least to lower order, see \cite{BegMa:qua}.

\section{Quasialgebras by cochain twist}

The constructions in the paper come out of quantum group theory (i.e. we use the language of Hopf algebras) but we apply them to classical (not quantum) enveloping algebras and finite group algebras.
Thus, let $H$ be a Hopf algebra with coproduct $\Delta:H\to H\tens H$, counit  $\eps:H\to \C$ and antipode $S:H\to H$, see \cite{Ma:book}. Let $F\in H\tens H$ be a cochain, i.e. $F$ is invertible and $(\eps\tens \id)F=1=(\id\tens\eps)F$. Associated to $F$ is its nonAbelian cohomology coboundary
\[ \Phi=\del F=F_{23}((\id\tens\Delta)F)( (\Delta\tens\id)F^{-1}) F_{12}^{-1}\]
where $F_{23}=1\tens F\in H^{\tens 3}$, etc. By construction $\Phi$, called the `associator', is a 3-cocycle in the required sense. These data go back to V.G. Drinfeld and it is known that $H^F$ defined
by the same algebra as $H$ and with coproduct $\Delta_F=F(\Delta\ )F^{-1}$ and suitable $S_F$ gives a quasi-Hopf algebra \cite{Dri:qua}.

Now let $A$ be an $H$-covariant associative algebra. The cochain-twisted quasialgebra $A_F$ is defined as the same vector space as $A$ but with a new product
\[ a\bullet b= \cdot(F^{-1}\la(a\tens b))\]
where $\la$ denotes the action of each copy of $H$. The new $A_F$ is nonassociative but obeys
\[( a\bullet b)\bullet c=\bullet (\id\tens (\  \bullet\  ))(\Phi\la (a\tens b\tens c))\]
for all $a,b,c$, and is covariant under $H^F$.

Moreover, when $\Omega(A)$ is an algebra of differential forms on
$A$ that is $H$-covariant, then $\Omega(A_F)=\Omega(A)_F$ defines
for us the wedge product algebra of differential forms on $A_F$,
covariant under $H^F$ and again potentially non-associative\cite{AkrMa:bra}. Note
that $\extd$ is not deformed and assumed to be commute with the
action of $H$, hence
\[ a\bullet \extd b=F^{-(1)}\la a \extd (F^{-(2)} b),
\quad \extd a\bullet b=\extd( F^{-(1)}\la a) F^{-(2)}\la b, \quad
\extd a\bullet \extd b= (\extd F^{-(1)}\la a)\wedge \extd
(F^{-(2)}\la b)\] for the deformed wedge product in terms of the
undeformed one, where $F^{-1}=F^{-(1)}\tens F^{-(2)}$ (summation
understood) is a notation.

The two examples that will be fully computed in the paper are of the general types which we now describe. Note that we work over $\C$ for convenience and because in physical examples there are  further unitarity restrictions (otherwise, the general constructions work over any field, though one should avoid certain characteristics in the examples). Also, we use the $H$-module version of the cochain twist theory as above because actions are more familiar to physicists; there is a parallel and in many ways better version of the theory with $H$ coacting on the algebra. 

\subsection{Quasi-$\R^n$}

Let $H=U(\R^n)$, with Hopf algebra structure
\[ \Delta \del^i=\del^i\tens 1+1\tens \del^i,\quad \eps \del^i=0,\quad S\del^i=-\del^i.\] Here $\R^n$ acts on $\R^n$ by translation and hence on its coordinate algebra $A=\C[\R^n]$ by differentiation operators $\del=\{\del^i\}$ and we think of the latter quite concretely as generating $U(\R^n)$. Let $F$ be a nowhere vanishing function of two vector coordinates (i.e. a function on $\R^{2n}$)  with value 1 when either argument is zero.  We consider  $F\in H\tens H$ (or in some completion of this space if $F$ is not a polynomial) as a cochain. Because $H$ is commutative, $H^F=H$ as an algebra and as a coalgebra, but is still regarded with
\[ \Phi(\del_1,\del_2,\del_3)=\frac{F(\del_2,\del_3)F(\del_1,\del_2+\del_3)}{F(\del_1+\del_2,\del_3)F(\del_1,\del_2)}\]
as a quasi-Hopf algebra. Here $\del_1=\del\tens 1\tens 1$, $\del_2=1\tens \del\tens 1$, $\del_3=1\tens 1\tens\del$ in $H^{\tens 3}$ so $\Phi$ is a function of these $3n$ variables.

 Then $A_F$ has a new product
 \[ a\bullet b=\cdot F^{-1}(\del_1,\del_2)a\tens b\]
 where $a(x)$, $b(x)$ are acted upon by $\del_1,\del_2$ respectively and then the result multiplied. Quasi-associativity will take the form above, as
 \[ (a\bullet b)\bullet c= \bullet (\id\tens (\  \bullet\  ))\Phi(\del_1,\del_2,\del_3)(a\tens b\tens c)\]
 where $\del_1$ means $\del$ acting on $a$, $\del_2$ means $\del$ acting on $b$, $\del_3$ means $\del$ acting on $c$, and products are in $A_F$. Recall that $\del$ itself is a vector, namely the momentum vector operator generating translations in $\R^n$.

Of interest in physics seems to be the following special case. Let $\square=\sum \del^i\tens{ \del^j} \eta_{ij}=\del_1\cdot \del_2$
taken with the Euclidean metric say (or any other fixed tensor $\eta$ on $\R^n$ in place of the dot product).  Let $f$ be any nowhere vanishing function in {\em one} variable and take
\[ F(\del_1,\del_2)=f(\square),\quad \Phi(\del_1,\del_2,\del_3)=\frac{f(\square_{23})f(\square_{12}+\square_{13})}{f(\square_{13}+\square_{23})f(\square_{12})}\]
where $\square_{13}=\del_1\cdot\del_3$ is $\square$ embedded in the first and third tensor positions, etc.

If $f$ is an exponential then $\Phi=1$ and $A_F$ is associative. For example, of $\eta_{ij}$ is antisymmetric one has the usual Moyal product for the Heisenberg-Weyl algebra or so-called noncommutative $\R^n$ used for example by  Seiberg and Witten for the effective description of the ends of open strings on 2-branes. At the other extreme would be $\eta_{ij}$ the Euclidean metric in which case the algebra remains commutative and associative. In general if $F$ remains symmetric but $f$ is no longer an exponential then the algebra $A_F$ will be commutative but not associative. This covers the example in \cite{Ram:gau} where
\[ F(\del_1,\del_2)=(1+ \frac{\lambda}{m} \square)^{-m}\]
which becomes approximately an exponential $\exp(\lambda \square)$ as $m\to \infty$. Here $\lambda$ is the deformation parameter which is taken with value $m^{-1}$ in \cite{Ram:gau}, but one can also keep these parameters $\lambda, m$ independent. We have
\[ \Phi(\del_1,\del_2,\del_3)=(1+\frac{ \frac{\lambda^2}{m^2}\square_{13}(\square_{12}-\square_{23})}{ 1+\frac{\lambda}{m}(\square_{12}+\square_{13}+\square_{23})+\frac{\lambda^2}{m^2}\square_{23}(\square_{12}+\square_{13})})^m\]

Another interesting family of commutative but nonassociative quasi-$\R^n$ is with
\[ F(\del_1,\del_2)=e^{-\frac{\lambda}{2}\square^2},\quad \Phi=e^{-\lambda \square_{13}(\square_{12}-\square_{23})}=e^{-\lambda \eta_{ij}\eta_{kl}(\del^i\del^k\tens \del^l\tens\del^j-\del^i\tens\del^k\tens\del^j \del^l)}\]
when we unpack our compact notation (summation convention understood).

A third variant is with $H=U(\R^n\lcross\R)$ where an extra 'dilation' generator $D$ is added. Its relations, coproduct and action on coordinates are
\[ [D,\del^i]=-\del^i,\quad \Delta D=D\tens 1+1\tens D,\quad D\la x_i=x_i\]
(so that $D$ has action $p$ on a monomial of total degree $p$). In this way $A=\C[\R^n]$
is again covariant under this extended $H$. One can now have more interesting cochains, for example
\[ F=e^{-\lambda\square- v(D\tens D)}\]
for a `potential function' $v$. If $v=0$ we have $\Phi=1$ as explained above.
In general is tempting to think of the introduction of non-bilinears in the exponent of $F$ as a way to encode interactions as non-associativity. The passage from the free theory to the interacting theory would then be a matter of a  cochain twist by the interaction\cite{AkrMa:bra}. This last example is in that spirit.

Clearly a great many models along the above lines are equally
possible, as any cochain $F$ is allowed in our framework.

\subsection{Quasi-$\Z_2^n$}

Here we take $H=\C(G)$, the functions on a finite group. This has basis of delta-functions $\{\delta_a\}$ labelled by
$a\in G$ and coproduct $\Delta \delta_a=\sum_{bc=a} \delta_b\tens \delta_c$, counit $\eps\delta_a=\delta_{a,e}$ and antipode $S\delta_a=\delta){a^{-1}}$. Here $e$ is the group identity.
We take $A=\C G$ the group algebra of $G$. This has basis $\{e_a\}$ labelled again by group elements. The product is just the product of $G$, so $e_ae_b=e_{ab}$.  This is covariant under $\C(G)$ with action
\[ \delta_a\la  e_b=\delta_{a,b}e_b.\]

A cochain on $H$ is a suitable $F\in H\tens H$ i.e. a nowhere
vanishing 2-argument function $F(a,b)$ on the group with value $1$
when either argument is the group identity $e$. Then
\[ \Phi(a,b,c)=\frac{F(b,c)F(a,bc)}{F(ab,c)F(a,b)}\]
is the usual group-cohomology coboundary of $F$ and is a group
3-cocycle. Then $H^F$ is the same algebra and coalgebra as $H$ but
is viewed as a quasi-Hopf algebra with this $\Phi$. Finally, the
canonical example of a quasi-algebra here is the twisted group
algebra  $A_F$ with the new product
\[ e_a\bullet e_b= F^{-1}(a,b) e_{ab}\]

An example is $G=\Z_2^3$ which we write additively as $3$-vectors $\vec{a}$ with entries in $\Z_2$.
We take
\[ F(\vec a,\vec b)=(-1)^{\vec a{}^T \begin{pmatrix}1 & 1 & 1 \cr  0 & 1 & 1\cr 0 & 0 & 1 \end{pmatrix} \vec b+ a_1b_2b_3+b_1a_2b_3+b_1b_2a_3},\quad \Phi(\vec a,\vec b,\vec c)=(-1)^{\vec a\cdot( \vec b\times\vec c)} \]
The new product \[ e_{\vec a}\bullet e_{\vec b}=F(\vec a,\vec
b)e_{\vec a+\vec b}\] is that of the octonions $\Bbb O$ as
explained in \cite{AlbMa:qua}. If we think of this in the same
spirit as the models above, we note that $\Phi$ comes from the
cubic `interaction term' in the exponent of $F$. Thus the
octonions are a cochain quantisation of the finite group $\Z_2^3$
as a quasi-algebra. Without the cubic interaction term one has the clifford algebra in
3 dimensions. Similarly $n=2$ gives the quaternions or (over $\C$) the algebra of $2\times 2$ matrices.

One can do the same for larger $\Z_2^n$. For the same bilinear
form as the above one, one obtains the Clifford algebra as an
associative cochain quantisation of $\Z_2^n$, while further `interaction'
terms give higher Cayley-Dickson and other quasi-algebras of
interest, see \cite{AlbMa:qua, AlbMa:cli}. Many other examples could be of
interest, eg for $G=\Z_n$ see \cite{AlbMa:Zn}.

\section{Gauge theory in monoidal categories} 

With the above background the
main question we address in this paper is that of gauge theory on
nonassociative spaces. For the ones in Section 2.1 of interest in
string theory, a somewhat complex approach has been proposed in \cite{Ram:gau} 
whereas here we propose a simpler one. Briefly, geometry including
gauge theory can be done in any monoidal Abelian category
$\CC$\cite{Ma:dia}\cite{Ma:rem}. We explain this in Section~3.1 and
give a concrete algebraic setting in Sections 3.2 and~3.3, which are the 
new results of the section. 

Before doing this, let us explain the problem at the simplest level. If we have an associative
algebra with a differential calculus obeying the Leibniz rule, one can write 
down the simplest '$U(1)$-Yang-Mills' theory where a connection is a differential 1-form
$\alpha\in \Omega^1$, decreed to transform as 
\eqn{usualgauge}{\alpha\mapsto \gamma^{-1}\alpha \gamma+\gamma^{-1}\extd \gamma} for $\gamma$ any invertible element of the algebra. The fundamental lemma of gauge theory is that then
the curvature $F(\alpha)=\extd \alpha+\alpha\wedge\alpha$ transforms by conjugation to $\gamma^{-1}F(\alpha)\gamma$. Note that
the non-linear term need not vanish in noncommutative geometry even in this simplest case.
The moduli space of flat connections up to gauge transformations is highly nontrivial even for the simplest commutative or noncommutative algebras \cite{MaRai:ele} and carries a lot of `homotopy' information. We will describe  it for  the functions on the cube $\Z_2^3$ in Section 5.3 under a further
unitarity restriction (in the $*$-algebra case one requires $\gamma^*=\gamma^{-1}$ i.e. unitary.)

Let us try this now when the algebra is nonassociative. The simplest part of the above lemma is that
$\alpha=\gamma^{-1}\extd \gamma$ should have zero curvature. Being careful about brackets, we have
\[ \extd(\gamma^{-1}\extd \gamma)+(\gamma^{-1}\extd \gamma)(\gamma^{-1}\extd \gamma)=(\extd \gamma^{-1})\extd \gamma+(\gamma^{-1}\extd \gamma)(\gamma^{-1}\extd \gamma)
=-((\gamma^{-1}\extd \gamma)\gamma^{-1})\extd \gamma+(\gamma^{-1}\extd \gamma)(\gamma^{-1}\extd \gamma)\]
which is nonzero precisely when $\gamma^{-1}\extd \gamma, \gamma^{-1},\extd \gamma$ fail to associate. The computation of $\extd \gamma^{-1}$ here is from $\extd (\gamma^{-1}\gamma)=0$ and the Leibniz rule, being careful about brackets.
This could work for some $\gamma$ in the algebra  but not for all invertible or unitary elements as in the associative case.   

\subsection{Diagrammatic gauge theory}

A monoidal category $\CC$ means a collection of objects with a tensor product
between any two objects and an associator natural isomorphism
$\Phi_{V,W,Z}:(V\tens W)\tens Z\to V\tens (W\tens Z)$ for any three objects, obeying the
usual properties, notably Mac Lane's pentagon identity. The latter says
that the two routes to rebracket 
\[ ((U\tens V)\tens W)\tens Z\to U\tens (V\tens (W\tens Z))\]
 are the same. In that case the
coherence theorem of Mac Lane says that all other bracketting
ambiguities are resolved, i.e. we can and should freely insert
associators $\Phi$ in order for expressions to make sense and
different ways to do that will give the same result. In that case
we can adopt a diagrammatic notation in which we omit brackets
entirely. We also denote $\tens$ by omission. We write maps
between objects (morphisms) as beads on a string flowing down from
one object to the other. We also require direct sums $\oplus$ to be defined and to be compatible in the usual way with $\tens$. Now,  because brackets are omitted, gauge theory {\em must} work at this
level because usual associative gauge theory works when expressed by the same diagrams.
In the nonassociative case, however, the translation of the diagrams back into algebra requires the insertion of the nontrivial associator $\Phi$ for rebrackettings. We recall here only the 'basic level' of gauge theory\cite{Ma:rem} in this diagrammatic form; there is a more geometrical theory with diagrammatic principal bundles etc.\cite{Ma:dia}. 

As an example an associative algebra $\CA$ in a monoidal category means an object $\CA$ with a product Y such that the two ways to feed the result of Y into another Y give the same. As a result we can depict the iterated product as a node with three lines coming in and one coming out (i.e. collapse the two equivalent tree graphs). We will use such a notation. A coalgebra $\CB$ is an object $\CB$ with a coproduct $\Delta:\CB\to\CB\tens\CB$ which we denote by an up-side-down Y and which `coassociates' similarly. The unit axiom for an algebra says that a 1 branching into a product can be `pruned' off. SImilarly a counit $\eps:\CB\to \und 1$ (the latter denoted by omission) is a branch emerging from a coproduct node and can be pruned. More details of 'algebra' in such diagrams are in \cite{Ma:book}. A coalgebra $\CB$ can 'coact' on an object $\CV$ and we use the up-side-down Y also to denote the coaction $\CV\to \CV\tens\CB$.

Similarly, a differential calculus $\Omega$ on $\CA$ means a graded algebra in the category with $\CA$ in degree zero, and $\extd$ a morphism (hence a node) increasing degree by 1, obeying a graded-Leibniz rule and $\extd^2=0$. All of this translates directly into (sums of) diagrams. One usually assumes that $\Omega$ is generated by $\CA$ and the 1-forms $\Omega^1$ but this is not necessary for the basic level of gauge theory that we describe here. We use Y also to denote products in this exterior algebra.

\begin{figure}\[\includegraphics{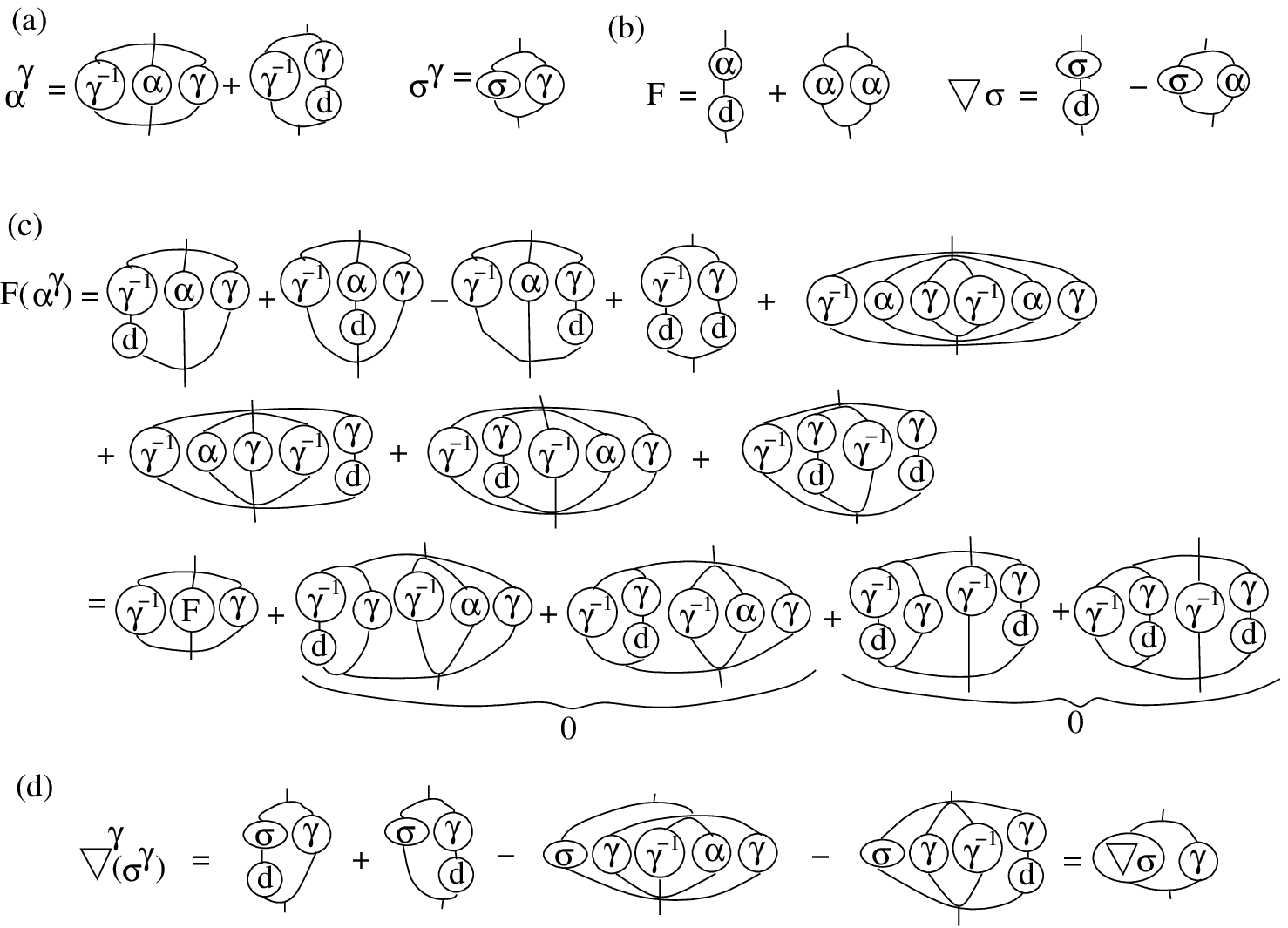}\]
\caption{Local gauge theory in a monoidal category: (a) gauge transform by $\gamma$ of gauge and matter fields (b) definition of curvature and covariant derivative and (c),(d) proof of covariance of $F,\nabla$}
\end{figure}

We are now ready to define matter fields as morphisms $\sigma: \CV\to \CA$. One can consider  that $\sigma$ has `values in $V^*$' (but it is more convenient to view it is a morphism). Similarly, a gauge field is a morphism $\alpha:\CB\to \Omega^1$ where $\CB$ is at least a coalgebra. Typically it might be a Hopf algebra in the category if this is braided, but such an assumption is again not needed for the basic level of gauge theory. One may think of $\alpha$ as a 1-form with values in the algebra $\CB^*$, i.e. we do general possibly non-Abelian gauge theory here, but again it is more convenient to view $\alpha$ as a morphism. Finally, a gauge transformation is a morphism $\gamma:\CB\to \CA$ with inverse $\gamma^{-1}$ in the sense
\[ \cdot(\gamma\tens\gamma^{-1})\Delta=\cdot(\gamma^{-1}\tens\gamma)\Delta=1\circ\eps\]
or in diagrams: if we split using the coproduct, apply $\gamma,\gamma^{-1}$ and close up with a product Y, this composition is the same either way as the counit map $\eps$ into nothing, followed by the unit map $1$ coming from nothing.   The action of such gauge transformations is shown in Figure~1(a). 
The basic covariant objects of interest namely the curvature and covariant derivative (the former is in a suitable sense the square of the latter) are shown in part (b) of the figure.

The fundamental lemmas of gauge theory are then shown in parts (c) and (d) of the figure; we check that $F(\alpha^\gamma)=F(\alpha)^\gamma$ and that $\nabla(\sigma^\gamma)=(\nabla\sigma)^\gamma$.  
In (c), we expand $\extd$ on the `conjugated' $\alpha$ using the Leibniz rule to obtain the first three terms. The next term is $\extd$ applied to `$\gamma^{-1}\extd\gamma$' again using the Leibniz rule, followed by $\extd^2=0$. The remaining four terms are an expansion of `$(\alpha^\gamma)^2$'.  Of the various terms, the 2nd and 5th (after cancelling $\gamma\gamma^{-1}$ to obtain a unit and counit and `pruning' these as explained above) give the transform of $F(\alpha)$ as required. The 1st (after inserting $\gamma\gamma^{-1}$ and 7th combine via Leibniz to give zero in view of $\extd(1)=0$. The 4th (inserting $\gamma\gamma^{-1}$) and 8th likewise give zero for the same reason. In (d) we compute $\nabla$ using the transformed quantities. The 2nd and 4th terms cancel (after cancelling $\gamma\gamma^{-1}$) and we identify the required result. 

 This establishes 'local gauge theory' at this diagrammatic level cf. \cite{Ma:rem} (where the focus was on the universal calculus, not assumed here). For principal bundles etc at this level see \cite{Ma:dia}. The latter contains explicit (associative) examples.

\subsection{Algebraic construction of nonassociative gauge theory by twisting}

The questions arise: how to obtain nonassociative examples of such a diagrammatic gauge theory and how does it look in explicit calculations? We will address the first in the remainder of the section, and the second in the remainder of the paper.

We do this by extending the cochain twisting theory in Section~2. Thus let $A$ be an algebra with calculus covariant under a background symmetry $H$ as in Section~2. Here $A$ could be functions on a classical manifold and $H$ the enveloping algebra of an ordinary Lie algebra, for example. Let  now $B$ be  an $H$-covariant coalgebra. It means that there is a coproduct $\Delta_B:B\to B\tens B$ which is an intertwiner for the aciton of $H$. Also a counit $\eps_B$. Suppose now that 
\[ \alpha:B\to \Omega^1(A),\quad \gamma:B\to A,\quad F(\alpha):B\to \Omega^2(A),\quad \sigma:V\to \Omega^1(A)\]
 and $\nabla$ are as in Section~3.1, i.e. a connection, gauge transformation etc. These form a gauge theory with the usual tensor product on the associative algebra $A$ as in Section 3.1. This theory can also be written without
diagrams by means of the `convolution product' $*$ of maps from a coalgebra to an algebra or module. Thus
\[ \alpha *\alpha= \wedge (\alpha\tens\alpha)\Delta_B,\quad F=\extd \alpha+ \alpha*\alpha,\quad \alpha^\gamma=\gamma^{-1}*\alpha*\gamma+\gamma^{-1}*\extd \gamma\]
and so forth. If $B=\C.1$ with $\Delta_B1=1\tens 1$, we have the simplest case of gauge theory mentioned in the preamble above. When $B=\C(G)$ the functions on a Lie group one has a general form of nonAbelian gauge theory. One can take here $H$ to be trivial, otherwise one has an equivariant gauge theory. 

Now let $F\in H\tens H$ be a cochain and define $B_F=B$ as a vector space but with deformed coproduct \[ \Delta_\bullet=F\la\Delta_B\]
 and unchanged $\eps_B$. Firstly, it can be seen that
$B_F$ is covariant under the twisted $H^F$. Indeed,
\[ h\la(F\la \Delta_B b)=F(\Delta h)F^{-1}\la(F\la\Delta_Bb)=F\la((\Delta h)\la\Delta_B b)=F\la\Delta_B (h\la b)\]
as the quasi-Hopf algebra $H^F$ acts on tensor products by its twisted coproduct $\Delta_F$ as explained in Section~2. Moreover, $B_F$ is a coalgebra but only in the monoidal category of $H^F$-modules, i.e. a 'quasi-coalgebra' in the sense:
\[\Phi_{B,B,B}(\Delta_\bullet\tens\id)\Delta_\bullet=(\id\tens\Delta_\bullet)\Delta_\bullet\]
as may be verified by direct computation. The theory is dual to that of twisting algebras so we omit the details. Similarly if $\Delta_V:V\to V\tens B$ is a coaction covariant under $H$, we define $V_F$ to be the same vector space but with deformed coaction $\Delta_{V\bullet}=F\la\Delta_V$, and can check that is is covariant under $H^F$ and a coaction of $B_F$ in the monoidal category.

We now claim that the same maps viewed as morphisms
\[ \alpha:B_F\to \Omega^1(A_F),\quad \gamma:B_F\to A_F,\quad F(\alpha):B_F\to \Omega^2(A_F),\quad \sigma:V_F\to A_F\]
form a gauge theory in the monoidal category of $H^F$-covariant objects, i.e. are an example of the constructions in Section~3.1 and enjoy the same relationships as before twisting. For example,
if we compute $\alpha *_\bullet\alpha$ where the subscript means in the deformed nonassociative theory,
\[ \alpha*_\bullet\alpha=\bullet(\alpha\tens\alpha)\Delta_\bullet=\bullet(F^{-1}\la (\alpha\tens\alpha)F\la\Delta_B=\wedge(\alpha\tens\alpha)\Delta_B=\alpha*\alpha\]
because each $\alpha$ is an intertwiner i.e. covariant under the action of $H$. We use $\bullet$ for the deformed product in the exterior algebra including wedge products. Similarly for all other expressions. In other words {\em the twisted non-associative theory is fully equivalent to the original associative one}. This is an important requirement from a deformation-theoretic point of view; if one thinks of the twisting as quantisation, this is an extension of the correspondence principle from classical gauge theory to gauge theory on the quantum (possibly nonassociative) space. 

On the other hand, computed entirely in the nonassociative deformed category, the gauge theory appears quite different. Remembering that the products are quasi-associative, we must fix brackettings when translating the diagrams into algebra and we do so by a convention to bracket by default to the left, inserting associators $\Phi$ according to Mac Lane's coherence theorem whenever a different bracketting is needed. Thus for example,
\[ \alpha^\gamma_\bullet=((\ \bullet\ )\tens\bullet\ )((\gamma^{-1}\tens\alpha)\tens\gamma)(\Delta_\bullet\tens\id)\Delta_\bullet+\bullet(\gamma^{-1}\tens\gamma)\Delta_\bullet\]
as a morphism $B_F\to \Omega^1(A_F)$. Provided one inserts $\Phi$ as specified (and where there is more than one way to do it one has the same result for any choice), the diagrammatic proof in Section~3.1 becomes an algebraic proof that 
\[ F_\bullet(\alpha)=\extd \alpha+\bullet(\alpha\tens\alpha)\Delta_\bullet\]
obeys
\[ F_\bullet(\alpha^\gamma_\bullet)=((\ \bullet\ )\tens\bullet\ )((\gamma^{-1}\tens F_\bullet)\tens\gamma)(\Delta_\bullet\tens\id)\Delta_\bullet\]
i.e. the fundamental lemma of (nonasociative) gauge theory. When there are matter fields we have similarly
\[ \sigma_\bullet^\gamma(v)=\bullet(\sigma\tens\gamma)\Delta_\bullet,\quad \nabla_\bullet\sigma(v)=\extd\sigma(v)-\bullet(\sigma\tens\alpha)\Delta_\bullet,\quad\nabla_\bullet^\gamma(\sigma_\bullet^\gamma)=(\nabla_\bullet\sigma)^\gamma.\]

\subsection{Canonical example equivalent to $U(1)$-Yang-Mills}

Finally, let us give a canonical example of an equivariant gauge theory and its twisting, that requires only the data for a cochain quantisation as in Section~2, i.e. there is a canonical choice of $B$. 

Thus, let $H$ be a Hopf algebra and $A$ and algebra with calculus which is $H$-covariant. We then set $B=H$ as a coalgebra, $\Delta_B=\Delta$ (the coproduct of $H$, ignoring the algebra structure of $H$). This automatically covariant under the action of $H$ on $B$ by left-multiplication: \[ h\la\Delta_B (b)=\Delta_H(h)\la\Delta_B(b)=\Delta_H(h)\Delta_H(b)=\Delta_H(hb)=\Delta_B(h\la b).\]
 
On the other hand, since every element of $B$ is obtained by acting by $H$ on 1, and since $\alpha,\gamma,F$ etc are morphisms, they are fully determined by their values on $1$, i.e. by
\[ \alpha(1)\in \Omega^1(A),\quad \gamma(1)\in A,\quad F(\alpha)(1)\in \Omega^2(A).\]
Here $\alpha(1),\gamma(1)$ etc. are chosen freely and form a usual gauge theory of the simplest $U(1)$-Yang-Mills type described in the preamble on any algebra. This is because $\Delta(1)=1\tens 1$ so all the coproducts in Figure~1 disappear when specialised to acting on  1, so
\[ \alpha^\gamma(1)=\gamma^{-1}(1)\alpha(1)\gamma(1)+\gamma^{-1}(1)\extd\gamma(1),\quad F(1)=\extd\alpha(1)+\alpha(1)\wedge\alpha(1)\]
etc. Our construction 'amplifies' this standard $U(1)$-Yang-Mills gauge theory on an algebra to an  $H$-equivariant one for any $H$  by $\alpha(b)=\alpha(b\la 1)=b\la\alpha(1)$ and $\gamma(b)=\gamma(b\la 1)=b\la\gamma(1)$.

For matter fields, the requirement that the coaction: $V\to V\tens B$ is a morphism makes $V$ into some form of 'Hopf module', i.e. a vector space on which $H$ both acts and coacts in a suitably compatible manner, namely here
\[ \Delta_V(h\la v)=(\Delta h)\la \Delta_V(v).\]
Hopf modules are fully determined by their space 
\[ V^H=\{v\in V\ |\ \Delta_V(v)=v\tens 1\}\]
of elements invariant under the coaction. The Hopf module-lemma ensures that these invariant elements $v\in V^H$ generate all of $V$ through the action. Note that this is usually done for action and coaction in the same side but with care works also in our case where the action is a left one and the coaction a right one. Indeed, we have
\[ H\tens V^H\to V,\quad h\tens v\to h\la v,\quad V\to H\tens V^H,\quad v\mapsto v^{(2)}{}_{(2)}\tens S^{-1}v^{(2)}{}_{(1)}\la v^{(1)}\]
where the antipode $S$ of the Hopf algebra is assumed to be invertible and where $\Delta_V(v)\equiv v^{(1)}\tens v^{(2)}$ and $\Delta h\equiv h_{(1)}\tens h_{(2)}$ are standard Hopf algebra notations. It is straightforward to see then that these two maps are mutually inverse, so $V\isom H\tens V^H$ and that the second map indeed lands in $H\tens V^H$ (this is not obvious but can be checked using routine Hopf algebra methods). Conversely, given any vector space $W$ we can define $V=H\tens W$ with action and coaction of $H$
\[h\la (g\tens w)=hg\tens v,\quad \Delta_V(h\tens w)=(h_{(1)}\tens w)\tens h_{(2)},\quad h,g\in H,\ w\in W\]
and check that  $W=V^H$; the above tells us that any crossed module $V$ is equivalent to one of this standard type. In short, the input data for matter fields in the theory boils down to choosing a vector space.

Moreover,  since $\sigma:V\to A$ is assumed to be $H$-covariant, it us fully determined by its values on this vector space $V^H$, since
$\sigma(\sum_p h_p\la v_p)=\sum_p h_p\la\sigma(v_p)$ for any basis $\{v_p\}$ of $V^H$. So the gauge theory above is equivalent to specifying a map $\sigma:V^H\to A$  or a multiplet of matter fields $\sigma(v_p)$ if we fix a basis of $V^H$. Thus our theory becomes equivalent to usual $U(1)$-theory with a multitplet of matter fields. Indeed,  $\sigma(v_p)\in A$ obeys
\[ \sigma^\gamma(v_p)=\sigma(v_p)\gamma(1),\quad (\nabla\sigma)(v_p)=\extd\sigma(v_p)-\sigma(v_p)\alpha(1)\]
as would be expected for $U(1)$ fields.

We now ready simply to twist this theory using the method in Section~3.2. $B_F$ now has 'deformed coproduct'  '$\Delta_\bullet=F\Delta$. A gauge field is again determined by $\alpha(1)$ but $\Delta_\bullet(1)=F\in H\tens H$ so 
\[F_\bullet (\alpha)(1)=\extd\alpha(1)+ (\alpha*_\bullet\alpha)(1)=\extd\alpha(1)+\bullet (\alpha\tens\alpha)(F) \]
\[ \alpha_\bullet^\gamma(1)=((\ \bullet\ )\bullet\ )((\gamma^{-1}\tens\alpha)\tens\gamma)((\Delta \tens\id)F)+\bullet(\gamma^{-1}\tens\extd\gamma)(F)\]
in terms of the deformed bullet product on $\Omega(A_F)$. As above, our convention is to read the diagrams with brackets accumulating to the left, with $\Phi$ to be inserted as needed for any other bracketting that may be required. The expressions above will be equal as linear maps  to $F(\alpha(1)), \alpha^\gamma(1)$, etc. as explained in Section~3.2, so the deformed theory is in correspondence with the original theory before twisting, but is well-formed in its own right. 

Finally,  if we have matter fields and elements $v_p$ that are invariant under the coaction, then the deformed coaction and hence gauge transform of matter fields is 
\[ \Delta_{V\bullet}(v_p)=F^{(1)}\la v_p\tens F^{(2)},\quad \sigma_\bullet^\gamma(v)=\sigma(F^{(1)}\la v_p)\gamma(F^{(2)}),\quad F\equiv F^{(1)}\tens F^{(2)}.\]
 Here we see that as with the gauge fields above, it is the entire 'amplified' theory that twists into a nonassociative one. It remains, however, equivalent to the $U(1)$-gauge theory with matter.

\section{Differentials and gauge theory on fuzzy $\R^n$}

In this section we illustrate the above formalism on the example of quasi-$\R^n$.
To be concrete, we focus calculations on the main example where
$f=(1+{\lambda\square\over m})^{-m}$ in Section~2.1, but the same methods
apply for the other versions of quasi-$\R^n$. We start with the algebra and
differentials in more detail, and then turn to the gauge theory. 

\subsection{Algebra and differentials on fuzzy $\R^n$} 

From Section~2.1, we have
\[ a\bullet b=\sum_{r=0}^m \left({m\atop r}\right) ({\lambda\over
m})^r(\del^{i_1}\cdots\del^{i_r}a) (\del_{i_1}\cdots\del_{i_r}
b),\] where we use $\eta_{ij}$ to lower indices. We call this
algebra $\R^n_{m,\lambda}$; the case in string theory is with
$\lambda={1\over m}$. For example, with the usual coordinates
$x_\mu$ of $\R^n$ we have the bullet product
\[ x_\mu\bullet x_\nu=x_\mu x_\nu + \lambda \delta_{\mu,\nu}\]
which is a simplified version \cite{Ram:gau} of higher-dimensional
fuzzy spheres that arise from the truncated matrix product in
certain string matrix models.  Here $m$ is a truncation order and the algebra becomes
associative as $m\to\infty$.

For our purposes we also need a differential calculus and we use
the same $F$ built from $f$ but now with the $\del^i$ acting by
Lie derivative on differential forms. Then the usual
$\Omega(\R^n)$ deforms to a (nonassociative)
$\Omega(\R^n_{m,\lambda})$. Notice that Lie derivative commutes
with exterior $\extd$, so the classical differential calculus is
indeed covariant as required. Then
\[ a\bullet\extd b=\sum_{r=0}^m \left({m\atop r}\right)({\lambda\over
m})^r(\del^{i_1}\cdots\del^{i_r}a) \extd\
\del_{i_1}\cdots\del_{i_r} b\]
\[\extd a\bullet\extd b=\sum_{r=0}^m \left({m\atop r}\right)({\lambda\over
m})^r\extd (\del^{i_1}\cdots\del^{i_r}a) \extd\
\del_{i_1}\cdots\del_{i_r} b\]  for functions $a,b$. 
For
example,
\[ x_\mu\bullet \extd x_\nu= x_\mu\extd x_\nu,\quad x_\mu\bullet \extd
(x_\nu\bullet x_\rho)=x_\mu\bullet \extd (x_\mu x_\nu)=x_\mu\extd
(x_\nu x_\rho)+\lambda(\delta_{\mu,\nu}\extd
x_\rho+\delta_{\mu,\rho}\extd x_\nu)\] \[ \extd x_\mu\bullet\extd x_\nu=\extd x_\mu\wedge \extd x_\nu=-\extd x_\nu\wedge\extd x_\mu=-\extd x_\nu\bullet\extd x_\mu\]
and so forth. This deformed
`quasidifferential calculus' is the classical one at lowest order
and but differentials of functions of degree $p$ will be modified
by descendants of lower degree. Because $\extd 1=0$ the relations involving
 $\extd x_\mu$ are necessarily unchanged,
\[ \extd a= (\del^\mu a)\extd x_\mu=(\del^\mu a)\bullet \extd x_\mu,\quad a\bullet \extd x_\mu=a\extd x_\mu=(\extd x_\mu)a=\extd x_\mu\bullet a.\]

\subsection{Gauge fields on fuzzy $\R^4$}

We are now ready to construct gauge theory on the above fuzzy $\R^4$ using the general construction in Section~3.3. 

First of all, we recall that here  $H=U(\R^n)=\C[\del^1,\cdots,\del^n]$  has coproduct $\Delta\del^i=\del^i\tens 1+1\tens\del^i$ on the generators.  We take for $B$ the same coalgebra, but to avoid confusion we denote this second copy  $B=U(\R^n)=\C[f^1,\cdots,f^n]$ with polynomial generators $f^i$. As before, we use a fixed (say Euclidean) $\eta_{ij}$ to lower indices. A gauge field is a covariant map $\alpha:B\to \Omega^1(\R^n)$ so  it is first of all a collection of 1-forms $\alpha(1),\alpha(f^i),\alpha(f^if^j)$ etc. in  $\Omega^1(\R^n)$. However, that $\alpha$ is a morphism requires
\[ \alpha(f^i)=\CL_i(\alpha(1))=\del^i\alpha(1)^\mu\extd x_\mu,\quad\cdots,\quad  \alpha(f^{i_1}\cdots f^{i_p})=\CL_{i_1}\cdots\CL_{i_p}(\alpha(1))=\del^{i_1}\cdots\del^{i_p}\alpha(1)^\mu\extd x_\mu.\]
where $\CL_i$ denotes the Lie derivative by the vector field $\del^i$ acting here on 1-forms. This is just action by $\del^i$ on the components $\alpha(1)^\mu$ in the coordinate basis.  This is how $\alpha(b)$ is determined from $\alpha(1)\in\Omega^1(\R^n)$. Similarly
\[ \gamma(f^i)=\del^i\gamma(1),\quad \cdots,\quad \gamma(f^{i_1}\cdots f^{i_p})=\del^{i_1}\cdots\del^{i_p}\gamma(1)\]
and similarly for $\gamma^{-1}$. This inverse is defined by the `convolution product', which involves the coproduct above, so for example
\[ \gamma^{-1}(1)\gamma(1)=1,\quad \gamma^{-1}(f^i)\gamma(1)+\gamma(1)\gamma(f^i)=0\]
\[ \gamma^{-1}(f^if^j)\gamma(1)+\gamma^{-1}(f^i)\gamma(f^j)+\gamma^{-1}(f^j)\gamma(f^i)+\gamma^{-1}(1)\gamma(f^if^j)=0\]
etc., which agrees with $\gamma^{-1}(f^i)=\del^i\gamma^{-1}(1)$ etc., as required by covariance. Similarly, we know that $\alpha^\gamma(1)=\alpha(1)^{\gamma(1)}=\alpha(1)+\gamma^{-1}(1)\extd\gamma(1)$. For higher order we compute the convolution product as
\begin{eqnarray*} \alpha^{\gamma}(f^i)&=&\gamma^{-1}(f^i)\alpha(1)\gamma(1)+\gamma^{-1}(1)\alpha(f^i)\gamma(1)+\gamma^{-1}(1)\alpha(1)\gamma(f^k)+\gamma^{-1}(f^i)\extd\gamma(1)+\gamma^{-1}(1)\extd \gamma(f^i)\\
&=&\alpha(f^i)+\gamma^{-1}(1)\extd\gamma(f^i)-\gamma^{-2}(1)\gamma(f^i)\extd\gamma(1)=\alpha(f^i)+\extd(\gamma^{-1}(1)\gamma(f^i))=\CL_i(\alpha^\gamma(1))\end{eqnarray*}
as it should as all our constructions are covariant under $H$.  Likewise, we know that $F(\alpha)(1)=F(\alpha(1))=\extd\alpha(1)$. At next order we have
\[ F(\alpha)(f^i)=\extd \alpha(f^i)+\alpha(f^i)\wedge\alpha(1)+\alpha(1)\wedge\alpha(f^i)=\extd\alpha(f^i)=\CL_i (F(\alpha)(1))\]
as it should. Thus the higher $\alpha(f^i)$ etc.,  behave like further auxiliary classical $U(1)$- gauge fields but are in fact determined from the $\alpha(1)$ gauge theory.  This gives the flavour of the amplified theory and its equivalence with usual $U(1)$ theory on $\R^n$. 

Next we deform to the coproduct of $B_F$,
\[ \Delta_\bullet f^i=(1+{\lambda\over m}f^j\tens f_j)^{-m}( f^i\tens 1+1\tens f^i)=f^i\tens 1+1\tens f^i-\lambda f^if^j\tens f_j-\lambda f^j\tens f_j f^i+\cdots\]
The action of $H$ on $B$ is multiplication by $\del^i=f^i$. 

As explained in Section~3.3 a gauge field still means an $H$-covariant map determined by $\alpha(1)\in\Omega^1(\R^n_{m,\lambda})$, i.e. some differential form $\alpha(1)=\alpha^\mu \bullet\extd x_\mu$.
Its curvature from Section~3.3 is
\begin{eqnarray*}F_\bullet(\alpha)(1)&=&\extd\alpha(1)+
\bullet(\alpha\tens\alpha)(F)=
\extd\alpha(1)+F^{(1)}\la\alpha(1)\bullet F^{(2)}\la\alpha(1)\\
&=& \extd\alpha(1)+\sum_{r=0}^\infty \left({m+r-1\atop r}\right)(-{\lambda\over
m})^r(\del^{i_1}\cdots\del^{i_r} \alpha^\mu\bullet\extd x_\mu)\bullet
(\del_{i_1}\cdots\del_{i_r} \alpha^\nu\bullet\extd x_\nu).\end{eqnarray*}
We know from the equivalence with the classical gauge theory that this will in fact equal $\extd\alpha(1)$ but this is a non-trivial computation from the point of view of the nonassociative theory. Similarly, we have
\[(\Delta\tens\id)(F)=(1+{\lambda\over m}(\square_{13}+\square_{23}))^{-m}\]
and hence
\begin{eqnarray*}\alpha^\gamma_\bullet(1)&=&((\ \bullet\ )\ \bullet\ )(1+{\lambda\over m}(\square_{13}+\square_{23}))^{-m}(\gamma^{-1}(1)\tens\alpha(1)\tens\gamma(1))\\
&& +\sum_{r=0}^\infty \left({m+r-1\atop r}\right)(-{\lambda\over
m})^r \del^{i_1}\cdots\del^{i_r} \gamma^{-1}(1)\bullet\extd 
\del_{i_1}\cdots\del_{i_r} \gamma(1). \end{eqnarray*}
where the first term can again be expanded as a powerseries as we have done for the second term. The action of a $\del^i$ on $\alpha$ is understood here to be via the Lie derivative. The second term is 'pure gauge' and we know by the equivalence with the untwisted theory that it is equal to $\gamma^{-1}(1)\extd\gamma(1)$ and hence its curvature is zero, as promised. From the point of view of the nonassociative theory, however,  these are nontrivial powerseries in the $\bullet$ product. Matter fields if present can similarly be included according to the theory at the end of Section~3.3. 

\section{Octonions as a finite
quasigeometries and gauge theory}

Here we illustrate the formalism of Section 3 on the octonions viewed as
a nonassociative coordinate ring obtained by quantising the classical space
$\Z^3_2$. The first section makes this point of view precise and is a main result
of the paper. We then consider gauge theory on this space.

\subsection{Octonions as quantisation and their differentials} 

The `classical' algebra of functions in the form of the group algebra $A=\C \Z_2^3$ before
deformation is generated by commuting $u,v,w$ say with
$u^2=v^2=w^2=1$. A general basis element is
\[ e_{\vec a}=u^{a_1}v^{a_2}w^{a_3}.\]
The deformed product has relations
\[ u\bullet u= v\bullet v= w\bullet w=-1,\quad u\bullet
v=-v\bullet u,\quad u\bullet w=-uw=-w\bullet v,\quad v\bullet
w=-wv=-w\bullet v\] which is indeed the usual octonions if one
puts $i=u,j=v$ and $k=u\bullet v$. Here
\[ F(\vec a,\vec a)=(-1)^{a_1+a_2+a_3+a_1a_2+a_1a_3+a_2a_3+a_1a_2a_3}=\begin{cases}1 & if \mbox{ $\vec a=0$} \\ -1 & else\end{cases},\]
which ensures
that $k^2=-1$ as it should. Similarly one may check that
\[k\bullet i=(u\bullet v)\bullet u=-(v\bullet u)\bullet
u=-v\bullet(u\bullet u)=v=j\] and so forth. Note that
\[ e_{\vec a}\bullet(e_{\vec b}\bullet e_{\vec c})=
(e_{\vec a}\bullet e_{\vec b})\bullet e_{\vec c}\] whenever $\vec
a,\vec b,\vec c$ are linearly dependent over $\Z_2$. This
expresses the `alternativity' property of the octions in our
formulation.

Next, the `classical' differential calculus on $A$ is fixed as
follows. By Fourier transform $A=\C(\hat\Z_2^3)$ where $\hat
\Z_2^3$ is position space if the previous $\Z_2^3$ above was momentum
space. Each $\hat \Z_2$ of position space is a finite set of two points
and it has only one possible differential calculus, the universal
one. It is then natural to take the three copies commuting (direct
product calculus), giving
\[ \extd u u= -u \extd u, \quad \extd u v= v\extd u,\quad \extd u
w= w \extd u\] and cyclic rotations. The wedge product is then
fixed by the graded Leibniz rule as
\[ \extd u \extd u=0,\quad \extd u \extd v=-\extd v\extd u,\quad
\extd u\extd w=-\extd w\extd v\] and cyclic rotations of this.
Notice that the more important objects here are the left-invariant
closed 1-forms
\[ \tau_1=-{1\over 2} u^{-1}\extd u,\quad \tau_2=-{1\over 2} v^{-1}\extd v,\quad
\tau_3=- {1\over 2} w^{-1}\extd w\] and the geometrical picture is
that of a 3-torus with the circle $S^1$ approximated by $\Z_2$. Moreover, the calculus has noncommutative de Rahm cohomology generated
by these $\tau_i$, exactly as for a classical 3-torus. These $\tau_i$ anti-commute among themselves in the wedge product and 
\[ \tau_ie_{\vec a}=(-1)^{a_i}e_{\vec a}\tau_i, \quad \extd e_{\vec a}= -2 e_{\vec a} a_i \tau_i.\]
We see that there is only a small amount of noncommutativity in our `classical' calculus attributable to the discrete nature of the underlying space.

The geometric picture here is clearer after making the above Fourier transform explicit. Thus, let
\[ e_{\vec a}(x)=(-1)^{a_i x_i};\quad u=(-1)^{x_1},\quad v=(-1)^{x_2},\quad w=(-1)^{x_3}\]
be the plane waves, where $ x $ is a point in position space (a $\Z_2$-valued vector). The exterior derivative here is 
\[ \extd f=( \del^i f )\tau_i,\quad \del^ie_{\vec a}=-2 a_i e_{\vec a}\]
where $\del^i$ is  the finite-difference operator in the $i$-direction. We see that the differentials act by multiplication in momentum space.

In general, the Fourier transform of $f(x)$ is a function $f_{\vec a}$ on momentum space characterised by $f(x)=\sum_{\vec a}f_{\vec a}e_{\vec a}(x)$. The inverse is
\[ f_{\vec a}={1\over 8}\sum_{ x } f(x)e_{\vec a}(x).\] 
Now, we have given the $\bullet$ deformation of $\Z_2^3$ into the octonions in momentum space as multiplication by $F(\vec a,\vec b)$. Let $F(x,y)$ be the same function before Fourier transform. Then
\begin{eqnarray*}( f\bullet g)(x)&=&\sum_{\vec a,\vec b}f_{\vec a}\  g_{\vec b} \  F(\vec a,\vec b)e_{\vec a+\vec b}(x)={1\over 64}\sum_{y,z,\vec a,\vec b}f(y)g(z)e_{\vec a}(y)e_{\vec b}(z)e_{\vec a+\vec b}(x)F(\vec a,\vec b)\\
&=&{1\over 64}\sum_{y,z}\sum_{y,z,\vec a,\vec b}f(y)g(z)e_{\vec a}(x+y)e_{\vec b}(x+z) F(\vec a,\vec b)\\ &=& {1\over 64}\sum_{y,z} F(y,z)f(x+y)g(x+z).\end{eqnarray*}
 Here
\begin{eqnarray*} &&\!\! F(y,z)=\sum_{\vec a,\vec b}(-1)^{a_1(b_1+b_2+b_3)+a_2(b_2+b_3)+a_3b_3+b_1 a_2 a_3+a_1 b_2 a_3+ a_1 a_2 b_3+ a_i y_i+ b_i z_i}\\
&&=2\sum_{a_1,a_2,b_2, b_3}(-1)^{(z_1+a_2a_3+a_2)(b_2+b_3)+a_3b_3+(z_1+a_2a_3)(b_2a_3+a_2b_3)+ (z_1+a_2a_3+z_1)y_1+a_2y_2+a_3y_3+b_2z_2+b_2z_3 }\\
&&=2^2 \sum_{a_3,b_3}(-1)^{(z_2+(z_1+z_2)a_3)b_3+a_3b_3+(z_1+z_2+a_3z_1)(z_1+a_3)b_3+y_1(z_1+a_3z_2)+y_2(z_1+z_2+a_2z_1)+y_3a_3+b_3z_3} \end{eqnarray*}
where we do the $b_1$ summation which gives a constraint $a_1+a_2a_3+z_1=0$ which eliminates $a_1$; then  we do the $b_2$ summation to obtain a constraint $a_2+z_1+z_2+a_3z_1=0$ to eliminate $a_2$. We next do the $b_3$ summation to obtain a constraint $a_3+z_1+z_2+z_1z_2+z_3=0$, giving
\[ F(y,z)=8 (-1)^{y^T\begin{pmatrix}1&1&0\\ 0&1&0\\ 1&1&1\end{pmatrix}z+y_1z_2z_3+z_1y_2z_3+z_1z_2y_3}.\]
We see that {\em the cochain $F$ that defines the octonions has the remarkable property that up to a relabelling, it is its own Fourier transform}, i.e. $F(y,z)$ has just the same form in position space as $F(\vec a,\vec b)$ in momentum space after a rotation of the indices $1\to 2\to 3\to 1$.
Note that the factor 8 in $F(y,z)$ is an artefact  due to our use of $1/8$ on one side of each Fourier transform rather than a  symmetrical $1/\sqrt{8}$. 

Note also that $f(x+y)=(R_1^{y_1}R_2^{y_2}R_3^{y_3})f(x)$ where $R_i$ is translation in the $i$ direction. Since $\del^i=R_i-1$, we have $f(x+y)=((1+\del^1)^{y_1}(1+\del^2)^{y_2}(1+\del^3)^{y_3}f)(x)$ which expresses the above result as a finite `bidifferential' operator
\begin{eqnarray*}
f\bullet g&=&\cdot({1\over 8}\sum_{y,z} (-1)^{y_1(z_1+z_2)+y_2z_2+y_3(z_1+z_2+z_3)+y_1z_2z_3+z_1y_2z_3+z_1z_2y_3}\\
&&\quad (1+\del^1)^{y_1}(1+\del^2)^{y_2}(1+\del^3)^{y_3}\tens (1+\del^1)^{z_1}(1+\del^2)^{z_2}(1+\del^3)^{z_3} )(f\tens g)\\
&=& \cdot\big(1\tens 1 - {1\over 2} (\del^1\tens \del^1 + \del^2\tens\del^1 + \del^3\tens\del^1  + \del^2\tens\del^2 +\del^3\tens\del^2+\del^3\tens\del^3\\
&&\quad  + \del^1\del^2\tens\del^1 + \del^1\del^3\tens\del^1 + \del^2\del^3\tens\del^1+\del^2\del^3\tens\del^2+ \del^2\tens\del^1\del^2\\
&&\quad+\del^3\tens\del^1\del^3 + \del^3\tens\del^2\del^3 +\del^1\del^2\del^3\tens\del^1+ \del^2\del^3\tens\del^1\del^2)\\
&& - {1\over 4}( -\del^1\tens\del^2\del^3    +\del^2\tens\del^1\del^3 +\del^3\tens\del^1\del^2 -\del^1\tens\del^1\del^2\del^3 +\del^2\tens\del^1\del^2\del^3 +\del^3\tens\del^1\del^2\del^3   \\
&&\quad +\del^1\del^2\tens\del^1\del^2 +\del^1\del^2\tens\del^1\del^3-\del^1\del^2\tens\del^2\del^3 + \del^1\del^3\tens\del^1\del^3 + \del^2\del^3\tens\del^1\del^3 \\
&&\quad +\del^2\del^3\tens\del^2\del^3 +\del^1\del^2\del^3\tens\del^1\del^2+\del^1\del^2\del^3\tens\del^1\del^3+\del^2\del^3\tens\del^1\del^2\del^3)\\
&& -   {1\over 8} \del^1\del^2\del^3\tens\del^1\del^2\del^3\big)(f\tens g)\end{eqnarray*}
These results have been obtained with MATHEMATICA. This makes precise the sense in which, in finite geometry, the octonions are a `quantisation' of functions on $\Z_2^3$. 

For comparison, if we do the same for the cochain that defines clifford algebras as a simpler associative quantisation of $\Z_2^n$, we have
\[ F(a,b)=(-1)^{a_1(b_1+\cdots+b_n)+a_2(b_2+\cdots+b_n)+\cdots +a_nb_n}\]\[  F(y,z)=2^n(-1)^{(y_1+y_2)z_1+(y_2+y_3)z_2+\cdots +(y_{n-1}+y_n)z_{n-1}+y_nz_n}.\]
The derivation of the latter is rather simpler than the above; we compute
\begin{eqnarray*}F(y,z)&=& (-1)^{a_1(b_1+\cdots+b_n)+a_2(b_2+\cdots+b_n)+\cdots +a_nb_n+\sum_{i=1}^na_iy_i+\sum_{i=1}^nb_iz_i}\\
&=&2(-1)^{z_1(y_1+y_2)} (-1)^{a_2(b_2+\cdots+b_n)+a_3(b_3+\cdots+b_n)+\cdots +a_nb_n+\sum_{i=2}^na_iy_i+\sum_{i=2}^nb_iz_i}\end{eqnarray*}
where we do the $b_1$ integral to obtain the constraint $a_1+z_1=0$, and change variables $a_2+z_1\to a_2$ in the result. What we obtain is $F(y,z)$ for $\Z_2^{n-1}$ in the remaining variables. The above then follows by induction. The $\bullet$-product description of the Clifford algebra in $n$-dimensions is then given as a quantisation of $\Z_2^n$ by this $F(y,z)$ by a similar formula as above. For example, for $n=3$ we have
\begin{eqnarray*}f\bullet g&=&\cdot({1\over 8}\sum_{y,z} (-1)^{(y_1+y_2)z_1+(y_2+y_3)z_2+y_3z_3} \\
&&\quad (1+\del^1)^{y_1}(1+\del^2)^{y_2}(1+\del^3)^{y_3}\tens (1+\del^1)^{z_1}(1+\del^2)^{z_2}(1+\del^3)^{z_3} )(f\tens g)\\
&=&\cdot\big(1\tens 1 - {1\over 2} (\del^1\tens \del^1 +\del^2\tens\del^2 + \del^3\tens\del^3 +\del^3\tens\del^2 +\del^2\del^3\tens\del^2 +1\tens\del^1\del^2\\
&&\quad  +\del^1\tens\del^1\del^2 +\del^2\tens\del^1\del^2  + \del^3\tens\del^1\del^2  + \del^3\tens\del^2\del^3  + \del^2\del^3\tens\del^1\del^2)\\
&& - {1\over 4} (\del^3\tens\del^1\del^2\del^3+\del^1\del^2\tens\del^1\del^2+\del^1\del^3\tens\del^1\del^2-\del^1\del^3\tens\del^1\del^3+\del^2\del^3\tens\del^2\del^3    \\
&&\quad  +\del^2\del^3\tens\del^1\del^2\del^3+ \del^1\del^2\del^3\tens\del^1\del^2) -{1\over 8} \del^1\del^2\del^3\tens\del^1\del^2\del^3\big)(f\tens g)\end{eqnarray*}

Finally, we turn to the differential geometry of the octonions. As a cochain twist we have that the relations involving the left-invariant  forms $\tau_i$ are unchanged (because $F$ acts trivially on them). Hence
\[ \extd e_{\vec a}=(\del^ie_{\vec a})\tau_i=(\del^i e_{\vec a})\bullet \tau_i,\quad e_{\vec a}\bullet \tau_i=e_{\vec a}\tau_i=(-1)^{a_i}\tau_i e_{\vec a}=(-1)^{a_i}\tau_i\bullet e_{\vec a}\]
in this basis.  For a more algebraic picture within the octonions, let us also consider
\[ E_{\vec a}=(u^{a_1}\bullet v^{a_2})\bullet w^{a_3}=(-1)^{a_1 a_2 a_3}u^{a_1}\bullet (v^{a_2}\bullet w^{a_3})=(-1)^{a_1a_2+a_1a_3+a_2a_3} e_{\vec a}\]
after a short computation using $F$ and $\Phi$. Since the $e_{\vec a}$ have square 1 with their initial product, and from the form of $F(\vec a,\vec a)$ above, we know that $E_{\vec a}\bullet E_{\vec a}=e_{\vec a}\bullet e_{\vec a}=-1$ with the exception of $E_{0}=1$. So these are all `unit octonions'.  Moreover, from the above,
\[ \extd E_{\vec a}=-2E_{\vec a} \bullet a_i\tau_i,\quad E_{\vec a}\bullet \tau_i =(-1)^{a_i}\tau_i\bullet E_{\vec a}\]
in this basis. We can then deduce
\[ \tau_1={1\over 2}u^{-1}\bullet\extd u,\quad \tau_2={1\over 2}v^{-1}\bullet\extd v,\quad \tau_3={1\over 2}w^{-1}\bullet\extd w\]
where inverse is in the octonions or bullet product algebra and eventually that
\[ \extd u\bullet u=-u\bullet\extd u,\quad \extd u\bullet v=-v\bullet\extd u,\quad \extd u\bullet w=-w\bullet \extd u\]
\[ \extd u\bullet\extd u=0,\quad \extd u\bullet\extd v=\extd v\bullet\extd u,\quad \extd u\bullet \extd w=\extd w\bullet\extd u\]
and cyclic rotations of this. The latter are obtained by applying $\extd$ and the graded-Leibniz rule which still holds. One can also obtain these results by direct computation from the action of $F$ and the initial calculus on $\Z_2^3$ as in \cite{AkrMa:bra}.

\subsection{Gauge fields on the octonions}

We have a basis of $H$ given by the $\delta$-functions $\{\delta_a\}$ on momentum space, with coalgebra
\[ \Delta\delta_{\vec a}=\sum_{\vec b+\vec c=\vec a}\delta_{\vec b}\tens\delta_{\vec c}\]
Their action on $A=\C\Z_2^3=\C(\hat \Z_2^3)$ is
\[ \delta_{\vec a}\la f(x)=f_{\vec a}e_{\vec a}(x),\quad \delta_{\vec a}\la(fg)=\sum_{\vec b+\vec c=\vec a}(\delta_{\vec b}\la f)(\delta_{\vec c}\la g)\]
i.e. it projects out the corresponding term in the Fourier expansion and behaves as shown on products.  We use the same coalgebra
$B$ with the same basis element $\delta_{\vec a}$ denoted  $f^{\vec a}$ to avoid confusion and the same form of coproduct as above. The action of $H$ is by $\delta_{\vec a}f_{\vec b}=\delta_{\vec a,\vec b}f_{\vec b}$. A gauge field is then a covariant map $\alpha:B\to \Omega^1(\hat\Z_2^3)$, i.e. a collection of 1-forms
\[ \alpha(f_{\vec a})=\alpha(\delta_{\vec a}\la 1)=\delta_{\vec a}\la\alpha(1)=(\delta_a\la\alpha(1)^i)\tau_i\]
where the $H$ acts trivially on the $\tau_i$ as explained in Section~5.1. Thus the collection is fully determined from $\alpha(1)=\sum_{\vec a}\alpha(f_{\vec a})$. Similarly the collection 
\[ \gamma(f_{\vec a})=\delta_{\vec a}\la\gamma(1),\quad \gamma(1)=\sum_{\vec a}\gamma(f_{\vec a})\]
is determined from the point-wise invertible function $\gamma(1)(x)$ in $\C(\hat \Z_2^3)$. The inverse $\gamma^{-1}(1)=\gamma(1)^{-1}$. More generally
\[\sum_{\vec b+\vec c=\vec a} \gamma^{-1}(f_{\vec b})\gamma(f_{\vec c})=\delta_{\vec a,0}\]
which is consistent with $\gamma^{-1}(f_{\vec a})=\delta_{\vec a}\la \gamma^{-1}(1)$. A gauge transform 
of $\alpha(1)$ is as usual
\begin{eqnarray*} \alpha^\gamma(1)&=&\alpha(1)^{\gamma(1)}=\gamma^{-1}(1)\alpha(1)\gamma(1)+\gamma^{-1}(1)\extd\gamma(1)\\
&=&\gamma(1)^{-1}(R_i\gamma(1))\alpha(1)^i\tau_i+\gamma(1)^{-1}\del^i\gamma(1)\tau_i=\alpha(1)+(\gamma(1)^{-1}\del^i\gamma(1))\phi(1)^i\tau_i\end{eqnarray*}
where there is a sum over $i$ and $\phi(1)^i=\alpha(1)^i+1$. Note that unlike the fuzzy $\R^n$ case the initial $U(1)$ theory already has a nontrivial conjugation because functions do not commute with the $\tau_i$ basic 1-forms. The change of variables to $\alpha=1+\phi$ is quite useful (see the next section) and $\phi$ transform by conjugation. For other components we have
\begin{eqnarray*}\alpha^\gamma(f_{\vec a})&=&\sum_{\vec b+\vec c+\vec d=\vec a}\gamma^{-1}(f_{\vec b})\alpha(f_{\vec c})\gamma(f_{\vec d})+\sum_{\vec b+\vec c=\vec a}\gamma^{-1}(f_{\vec b})\extd\gamma(f_{\vec c})\\
&=& \sum_{\vec b+\vec c+\vec d=\vec a}\gamma^{-1}(f_{\vec b})R_i(\gamma(f_{\vec d}))\alpha(f_{\vec c})^i\tau_i+\delta_{\vec c,0}\gamma^{-1}(f_{\vec b})\del^i\gamma(f_{\vec d})\tau_i\\
&=&\alpha(f_{\vec a})+\sum_{\vec b+\vec c+\vec d=\vec a}(\gamma^{-1}(f_{\vec b})\del^i\gamma(f_{\vec d}))(\alpha(f_{\vec c})^i+\delta_{\vec c,0})\tau_i=\delta_a\la(\alpha^\gamma(1))\end{eqnarray*}
as it should. For the last step we identify $\delta_{\vec c,0}=\delta_{\vec c}\la 1$ as $1=e_0(x)$ and use the action of $\delta_{\vec a}$ on triple product along the lines of its action on a product explained above. Thus the theory looks like a collection of 1-forms with gauge-like transformation properties but determined consistently from the single theory for $\alpha(1)$.  Similarly, for the curvature we have
\begin{eqnarray*}F(\alpha)(1)&=&F(\alpha(1))=\extd\alpha(1)+\alpha(1)\wedge\alpha(1)=\sum_{i,j}(\del^i\alpha(1)^j+\alpha(1)^iR_i\alpha(1)^j)\tau_i\wedge \tau_j\\
&=&\sum_{ij}(\del^i\alpha(1)^j+\alpha(1)^i\del^i\alpha(1)^i+\alpha(1)^i\alpha(1)^j)\tau_i\wedge \tau_j=\sum_{ij}\phi(1)^i\del^i\phi(1)^j\tau_i\wedge \tau_j\end{eqnarray*}
where $\alpha(1)^i\alpha(1)^j\tau_i\wedge \tau_j=0$ as the $\tau_i$ anticommute. This is a standard form for the $U(1)$-Yang-Mills curvature on a discrete space in noncommutative geometry. The other components may similarly be computed  as
\[F(\alpha)(f_{\vec a})=\extd\alpha(f_{\vec a})+\sum_{\vec b+\vec c=\vec a}\alpha(f_{\vec b})\wedge\alpha(f_{\vec c})=\sum_{ij}\sum_{\vec b+\vec c=\vec a}\phi(f{\vec b})^i\del^i\phi(f_{\vec c})^j\tau_i\wedge\tau_j=\delta_{\vec a}\la F(\alpha(1))\]
as it should, where $\delta_{\vec a}$ acts on the coefficients of $\tau_i\wedge\tau_j$, i.e. the other components have a similar form but are determined by $F(\alpha(1))$. 

The above `amplification' of $\alpha(1)$ to a collection of gauge fields  can be made even more explicit by different basis $e_y\equiv\sum_{\vec a}e_y(\vec a)\delta_{\vec a}$ of $H$ where $e_y(\vec a)=e_{\vec a}(y)=(-1)^{y_ia_i}$. These elements have $\Delta e_y=e_y\tens e_y$ (this is the isomorphism $\C(\Z_2^3)\isom \C\hat \Z_2^3$). They act on functions by $(e_y\la f)(x)=f(x+y)$ and $\alpha(e_{\vec x})$ behave more explicitly like $\alpha(1)$, which is one of the collection via $e_0=1$.

We now turn to the twisted nonassociative theory. The coproduct of $B_F$ is
\[ \Delta_\bullet f_{\vec a}=\sum_{\vec b+\vec c=\vec a}F(\vec b,\vec c)f_{\vec b}\tens f_{\vec c},\quad \Delta_\bullet E_x={1\over 64}\sum_{y,z}F(y,z)E_{x+y}\tens E_{x+z};\quad E_y=\equiv\sum_{\vec a}e_y(\vec a)f_{\vec a}\]
where the cochain and its Fourier transform are (if we wan the octonions) as in Section~5.1. As explained in Section~3.3 a gauge field still means an $H$-covariant map determined by $\alpha(1)\in\Omega^1(\Bbb O)$, i.e. some differential form $\alpha(1)=\alpha^i\bullet\tau_i=(\phi^i+1)\bullet\tau_i$ (sum over $i$). The curvature according to Section~3.3, is
\[ F_\bullet(\alpha)(1)=\extd\alpha(1)+\bullet(\alpha\tens\alpha)(F)=\extd\alpha(1)+\sum_{\vec b,\vec c}F(\vec b,\vec c)(\delta_{\vec b}\la \alpha(1))\bullet (\delta_{\vec c}\la\alpha(1)).\]
Similarly, we have
\[ (\Delta\tens\id)(F)=\sum_{\vec a, \vec b,\vec c} F(\vec a+\vec b,\vec c)\delta_{\vec a}\tens\delta_{\vec b}\tens\delta_{\vec c}=\sum_{\vec a, \vec b,\vec c} F(\vec a,\vec c)F(\vec b,\vec c)\delta_{\vec a}\tens\delta_{\vec b}\tens\delta_{\vec c}=F_{13}F_{23}\]
for the particular form of $F$ for the octonions (which is linear in the exponent with respect to the first argument). Then
\[ \alpha_\bullet^\gamma(1)= \sum_{\vec a, \vec b,\vec c} F(\vec a,\vec c)F(\vec b,\vec c)((\delta_{\vec a}\la\gamma^{-1}(1)\bullet \delta_{\vec b}\la\alpha(1))\bullet\delta_{\vec c}\la\gamma(1))+
\sum_{\vec b,\vec c}F(\vec b,\vec c)\delta_{\vec b}\la\gamma^{-1}(1)\bullet \extd \delta_{\vec c}\la\gamma(1).\]
Matter fields can similarly be included from the general theory in Section 3.3.

Next, although our view of the octonions as a nonassociative quantisation of functions on the cube is the 'geometrical one', it remains very convenient to work with our original plane-wave basis $\{e_{\vec a}\}$ for calculations. Here $f\in H$ acts diagonally by as multiplication by $f(\vec a)$ on an element of degree $\vec a$. Here the degree is multiplicative and $\extd $ does not change degree, so for example $e_{\vec a}\extd e_{\vec b}$ has degree $|e_{\vec a}\extd e_{\vec b}|=\vec a+\vec b$. Similarly after deformation with the $E_{\vec a}$. Here for example $E_{\vec a}\bullet\extd E_{\vec b}=F(\vec a,\vec b)e_{\vec a}\extd e_{\vec b}$ etc. The hard part from this point of view is to find the inverse in the undeformed algebra of a general gauge transformation $\gamma=\sum_{\vec a}\gamma_{\vec a}e_{\vec a}$.  The answer is to construct $\gamma^{-1}=\sum_{\vec a}\gamma^{-1}_{\vec a}e_a$  by Fourier transform of the inversion operation:
\[  \gamma^{-1}_{\vec a}={1\over 8}\sum_x {e_{\vec a}\over\gamma(x)}={1\over 8}\sum_x {e_{\vec a}\over\sum_{\vec b}\gamma_{\vec b}e_{\vec b}(x)}
\]
where we require the $\gamma(x)$ (the sum in the denominator) to be non-zero for each $x$, i.e. all signed sums of the $\gamma_{\vec b}$ coefficients should be non-zero. Otherwise, since the action of $\delta_{\vec a}$ on any expression in the exterior algebra is to pick out the degree $\vec a$ part, we have more simply in this basis:
\[ F_\bullet(\alpha)=\extd\alpha+\sum F(|\alpha|,|\alpha'|)\alpha\bullet\alpha'\]
\[ \alpha_\bullet^\gamma= \sum F(|\gamma^{-1}|,|\gamma|)F(|\alpha|,|\gamma|) (\gamma^{-1}\bullet\alpha)\bullet\gamma+\sum F(|\gamma^{-1}|,|\gamma|)\gamma^{-1}
\bullet\extd\gamma\]
where the sum is over the different graded components of each object and to this end $\alpha'$ denotes a second independent copy of $\alpha$. Also, we omit writing that these are the gauge and other fields at 1, i.e. $\alpha\equiv\alpha(1)$ etc. Even though the amplification to the collection of fields is needed for the diagrammatic picture of Section~3.1, all formulae are by now referred back to their values on 1. If one wants to be more explicit and write the homogeneous degree components explicitly, we have
\[ F_\bullet(\alpha)=\extd\alpha+\sum_{\vec a,\vec b}F(\vec a,\vec b)\alpha_{\vec a}\bullet\alpha_{\vec b}\]
\[ \alpha^\gamma_\bullet=\sum F(\vec a,\vec c)F(\vec b,\vec c) (\gamma^{-1}_{\vec a}E_{\vec a}\bullet\alpha_{\vec b})\bullet\gamma_{\vec c}E_{\vec c}+\sum_{\vec a,\vec b} F(\vec a,\vec b)\gamma^{-1}_{\vec a}E_{\vec a}\extd\gamma_{\vec b}E_{\vec b}\]
where $\alpha=\sum_{\vec a}\alpha_{\vec a}$ is the decomposition into homogeneous components (this is a slightly different notation from the Fourier decomposition of $\gamma$ into components $\gamma_{\vec a}E_{\vec a}$). The fuzzy-$\R^4$ example in Section~4.2 can likewise be computed more simply in this 'momentum space' point of view. 

Finally, we demonstrate this gauge theory with an example of a completely explicit computation. Thus, let 
\[\gamma=\lambda u+\mu v,\quad \gamma^{-1}={1\over\lambda^2-\mu^2}(\lambda u-\mu v)\]
where $u=e_{(1,0,0)}$, $v=e_{(0,1,0)}$ are two of the octonion generators as explained in Section~5.1, and $\lambda\ne\pm\mu$. These are necessarily also inverse in the convolution-algebra $\gamma^{-1}*_\bullet\gamma=1$ as one may verify directly. Similarly, since 
$F(|u|,|u|)=F(|v|,|v|)=F(|u|,|v|)=-1$, $F(|v|,|u|)=1$, we have 
\[ \gamma^{-1}*_\bullet\extd\gamma={1\over\lambda^2-\mu^2}(-\lambda^2 u\bullet\extd u-\lambda\mu(u\bullet\extd v+v\bullet\extd u)+\mu^2v\bullet\extd v)\]
Let us check that the curvature of this pure gauge part is zero:
\[ \extd (\gamma^{-1}*_\bullet\extd\gamma=-{\lambda\mu\over\lambda^2-\mu^2}(\extd u\bullet\extd v+\extd v\bullet\extd u)=-2{\lambda\mu\over\lambda^2-\mu^2}\extd u\bullet\extd v\]
using the relations in the octonion calculus from Section~5.1. Meanwhile, when we square $\gamma^{-1}*_\bullet\gamma$ in the convolution product we must insert the factors $F(\vec a,\vec b)$ when multiplying components of degrees $\vec a,\vec b$ as explained above. Here $|u\bullet\extd v|=|v\bullet\extd u|=(1,1,0)$ while $u\bullet\extd u$ and $v\bullet\extd v$ have degree 0. Hence of the 16 terms only four come in with a - sign. Moveover, when we multiply out the 16 terms we can, in these particular expressions, associate, because the degree vectors for $u,v,u\bullet \extd v,v\bullet\extd u$ are linearly independent, so $\Phi$ for them is trivial. This results in all but four of the terms zero or cancelling pairwise. For example
\[ (u\bullet\extd v)\bullet(u\bullet\extd v)=u\bullet((\extd v\bullet u)\bullet \extd v)=-u\bullet(u\bullet(\extd v\bullet\extd v))=0\]
\begin{eqnarray*}(u\bullet\extd v)\bullet (v\bullet\extd u)&=&u\bullet((\extd v\bullet v)\bullet\extd u)=-u\bullet(v\bullet(\extd v\bullet\extd u))=-(u\bullet v)\bullet(\extd u\bullet\extd v)\\
&=&(v\bullet u)\bullet(\extd u\bullet\extd v)=-(v\bullet\extd u)\bullet(u\bullet\extd v)\end{eqnarray*}
using the relations from Section~5.1 (the last step is analogous the first sequence). What remains is
\begin{eqnarray*} (\gamma^{-1}*_\bullet\extd\gamma)*_\bullet (\gamma^{-1}*_\bullet\extd\gamma)&=& 
{\lambda\mu\over(\lambda^2-\mu^2)^2}(\lambda^2((u\bullet\extd u)\bullet(u\bullet\extd v)+(u\bullet\extd v)\bullet(u\bullet\extd u))\\
&&-\mu^2((v\bullet\extd v)\bullet(v\bullet\extd u)+(v\bullet\extd u)\bullet(v\bullet\extd v)))\\
&=&2{\lambda\mu\over\lambda^2-\mu^2}\extd u\bullet\extd v\end{eqnarray*}
by similar computations 
\[ (u\bullet\extd u)\bullet(u\bullet\extd v)=u\bullet((\extd u\bullet u)\bullet\extd v)=-u\bullet ((u\bullet\extd u)\bullet \extd v)=-(u\bullet u)\bullet(\extd u\bullet\extd v)=\extd u\bullet\extd v\]
etc., using the relations of the octonion calculus. Hence the curvature of this pure gauge part is zero as promised.

\subsection{Moduli of zero-curvature $U(1)$-Yang-Mills connections on $\Z_2^n$ and octonions}

By construction the above example of gauge theory on the octonions (not the only possible one, depending on the choice of gauge group coalgebra), is equivalent to that in the 'classical' object $\Z_2^3$. Maxwell theory on $\Z_2^n$ (but not Yang-Mills) has been covered in \cite{Ma:phy} and also
quantum Yang-Mills theory on $\Z_2^2$ but the analysis for classical $U(1)$-Yang-Mills and in particular the moduli space of zero curvature solutions has not to our knowledge been given even for $\Z_2^2$. We fill this gap now. As to be expected on a torus, this moduli space is nontrivial.

We use the `classical' calculus on $\Z_2^n$ as described in Section 5.1 before deformation to the octonions when $n=3$. The exterior algebra is generated by the plane-wave functions $e_{\vec a}(x)$ (now $a$ an $n$-vector) and $\tau_i$, $i=1,\cdots n$ as in Section 5.1 but now for general $n$. They anticommute among themselves, etc.

A $U(1)$-Yang-Mills gauge field means $\alpha=\alpha^i\tau_i\in \Omega^1(\Z_2^n)$ where the $\alpha^i(x)$ are the component functions. The curvature $F=\extd\alpha+\alpha\wedge\alpha$ is
\[ F=\sum_{i<j}F^{ij}\tau_i\wedge\tau_j,\quad F^{ij}=\del^i\alpha^j-\del^j\alpha^i+\alpha^iR_i\alpha^j-\alpha^jR_j\alpha^i.\]
We change variables to $\phi^i=1+\alpha^i$ or $\alpha=\phi-\vartheta$, where $\phi=\phi^i\tau_i$ and  $\vartheta=\sum_i\tau_i$ is a zero curvature `reference' connection that is closed, not exact and squares to zero. The moduli of flat connections contains at least this nonzero gauge-invariant element. Indeed, let $\gamma$ be a point-wise invertible function on position space. The gauge transformation of $\alpha$ and the expression for the curvature in terms of $\phi^i$ are:
\[( \phi^\gamma)^i={\gamma\over R_i\gamma}\phi^i,\quad F^{ij}=\rho_i\phi^j-\rho_j\phi^i\]
where $\rho_i\equiv \phi^iR_i$. Finally, to be physical, we fix unitarity conditons. As in \cite{MaRai:ele} we require the $\tau_i$ to be invariant under a $*$-operation extending the point-wise complex conjugation operation on position space. We then require the $\alpha$ to be hermitian, which translates in view of the commutation relations between the $\tau_i$ and functions to
\[ \bar{\phi^i}=R_i\phi^i,\quad \gamma=e^{\imath \xi},\quad(\phi^\gamma)^i=e^{-\imath\del^i\xi}\phi^i.\]
The middle equation is because if the reality of all the $\phi^i$ is preserved one may deduce that $\del^i(\bar\gamma\gamma)=0$ so $\gamma$ is without loss of generality pointwise unitary. We then put this into the transformation of $\phi^i$.

By a similar argument to the proof for $S_3$ in \cite{MaRai:ele} we have for all $i,j$ in the case of a zero-curvature solution:
\begin{eqnarray*}\rho_i\rho_i\phi^j&=&\phi^iR_i(\phi^i R_i\phi^j)=\lambda_i^2\phi^j\\
&=&\rho_i\rho_j\phi^i=\phi^iR_i(\phi^jR_j\phi^i)=\rho_i(\phi^j) R_jR_i\phi^i=\rho_j(\phi^i)R_i R_j\phi^i=\phi^jR_j(\lambda_i^2)\end{eqnarray*}
where 
\[ \phi^i=\lambda_ie^{\imath\theta_i},\quad \lambda_i^2=|\phi^i|^2=\phi^iR_i\phi^i\]
is a polar decomposition. We conclude that at each point $R_i\lambda_i=\lambda_i$ by the reality and 
at each point either  $\lambda_i=0$ or  $R_i\lambda_j=\lambda_j$ by the above computation, i.e.
\[ \lambda_i \del^i\lambda_j=0,\quad\forall i,j,\quad  \del^i\lambda_i=0,\quad\forall i\]
 These $\lambda_i$ are gauge-invariant and we now use them to analyse the possible solutions.

{\bf case 1: $\exists$ a point with all $\lambda_i\ne 0$ (constant maximal case)}. In this case each $\lambda_j$ will be unchanged moving in every direction to an adjacent point. Hence at each adjacent point they will all be nonzero. We conclude that all the $\lambda_i$ are constant functions.

Moreover, in this case the zero curvature equations become 
\[ e^{\imath\theta_i}e^{\imath R_i\theta_j}e^{-\imath R_j\theta_i}e^{-\imath\theta_j}\]
after cancelling $\lambda_i\lambda_j$ from both sides. If we think of $e^{\imath\theta_i(x)}$ to be a factor for parallel transport along the edge in direction $i$ from $x$, then this says that the holonomy around the plaquet with bottom left corner $x$ and edges in the $i$ and $j$ directions is zero. In this case, from such a solution we construct the following gauge transform: 
\[ \gamma(0)=1,\quad \gamma(x)=e^{\imath\sum_a\imath\theta.\extd a}\]
where we take any path $a$ from 0 to $x$ and multiply the parallel transports on the edges of the path. As in usual gauge theory, this transforms all the $\phi^i\to \lambda_i$, i.e. eliminates all the phases as gauge degrees. Hence the solutions up to gauge equivalence in this case are of the form
\[ \alpha=\lambda_i\tau_i-\vartheta,\quad \lambda_i\in \R_{>0},\quad i=1,\cdots, n.\]

{\bf case 2: $\exists$ a point with exactly one $\lambda_i= 0$ (split case).} In this case all $\lambda_j\ne 0$ for $j\ne i$, at the point in question. Therefore moving in all directions other than $i$, we have the same value for all $\lambda_j$, i.e. a constant maximal solution on the  subspace $\Z_2^{n-1}$. We also have the same value of $\lambda_i=0$ throughout this subspace. Moreover, moving in the $i$ direction from any point in the subspace keeps $\lambda_i=0$ (hence $\lambda_i\equiv 0$ everywhere) but says nothing about the values of any of the $\lambda_j$, $j\ne i$. Hence the solution is two independent copies of  ${n-1}$-dimensional solutions, in which the first copy is maximal by assumption and the second copy is unconstrained. 

For example, if the second copy is also maximal on the subspace, we have
\[ \alpha=x_i\lambda_j\tau_j+(1-x_i)\mu_j\tau_j-\vartheta,\quad \lambda_j,\mu_j\in \R_{>0},\quad j\ne i\]
and again any solution in this case is equivalent to something of this form (one may gauge away the phases in each $\Z_2^{n-1}$ space separately). The value of $\del^i\xi$ between the two copies could produce a gauge phase factor but this is irrelevant as  $\lambda_i\equiv 0$ everywhere.

{\bf case 2'} What arises naturally here is the weaker assumption just that there is some $i$ with $\lambda_i\equiv0$ throughout the space. In this case the solution   necessarily splits into independent solutions of any type of one dimension lower, of whatever type.  This is therefore covered by
induction.  Hence it remains only to classify the remaining cases under the assumption  that the solution is not split in any direction.

{\bf case 3: $\exists$ a point with exactly two $\lambda_i=\lambda_j=0$ and no splitting.} Here as before there is a $\Z_2^{n-2}$ subset containing the point with $\lambda_k\ne 0$ for all $k\ne i,j$ and $\lambda_i=\lambda_j=0$ throughout. Moreover, stepping in the $i$ direction carries over $\lambda_i=0$ to the entire adjacent quadrant, but none of the other information. Similarly stepping in the $j$ direction carries over $\lambda_j=0$ to that quadrant.   We then relate the quadrants by further analysis; see the example below.

One may proceed in this way to classify the cases with more and more assumed degeneracy. Among the solutions are those of the same form as the constant maximal case above but allowing any of the $\lambda_i=0$. These are multiply-split solutions and include each $\tau_i$ alone as a zero-curvature flat connection, as well as $\alpha=\vartheta$.

To be concrete we now offer a complete classification for $n=2$ and $n=3$ which demonstrates the method. The $n=3$ case is  in correspondence with solutions on the octonions by twisting as we have mentioned.

{\bf For $n=2$} we have two cases: (i)  the constant maximal solution is 
\[ \alpha=\lambda_1\tau_1+\lambda_2\tau_2-\vartheta,\quad \lambda_i\in \R_{>0}\]

(ii) we  have a splitting $\lambda_1\equiv 0$,  with  $\lambda_2$ unconstrained other than being constant in the $2$-direction, i.e.
\[ \alpha=\lambda\tau_2-\vartheta,\quad\del^2\lambda=0\]
where $\lambda$ is a function just in the $x_1$ variable and up to gauge equivalence can be taken real and non-negative in its values. Similarly for a splitting $\lambda_2\equiv 0$:
 \[ \alpha=\lambda \tau_1-\vartheta,\quad \del^1\lambda=0\]
for a real and non-negative function $\lambda$ of $x_2$ alone.  
 
 {\bf For $n=3$} we have three cases: (i) the constant maximal solution
 \[ \alpha=\lambda_1\tau_1+\lambda_2\tau_2+\lambda_3\tau_3-\vartheta,\quad\lambda_i\in\R_{>0}\]
 
 (ii) we have a splitting $\lambda_1=0$:
 \[ \alpha=x_1\phi+(1-x_1)\psi-\vartheta \]
 where $\phi,\psi$ correspond to two independent solutions on the $\Z_2^2$ subsets (faces) with $x_1=0$ and $x_1=1$ respectively. Up to gauge transformation they can be taken real and positive, i.e. without phases.
 Similarly in the other two directions.  

 (iii) we suppose that there does not exist a splitting, but there does exist a point with, say, $\lambda_1=\lambda_2=0$ and $\lambda_3=\nu\ne 0$. To be concrete let this point be $A$ the origin
  in the standard cube shown in Figure~2. These are also the values at $H$ by the above argument; the equal value of $\lambda_3$ is shown in part (a) of the figure by labelling the arrowed edge, and such a nonzero edge 'transports' the other values from  $A$ to $H$ by the arguments above. We also see that $\lambda_1=0$ at $B$ and $G$, while $\lambda_2=0$ at $D$ and $E$, by the reality condition.

   \begin{figure} \[ \includegraphics{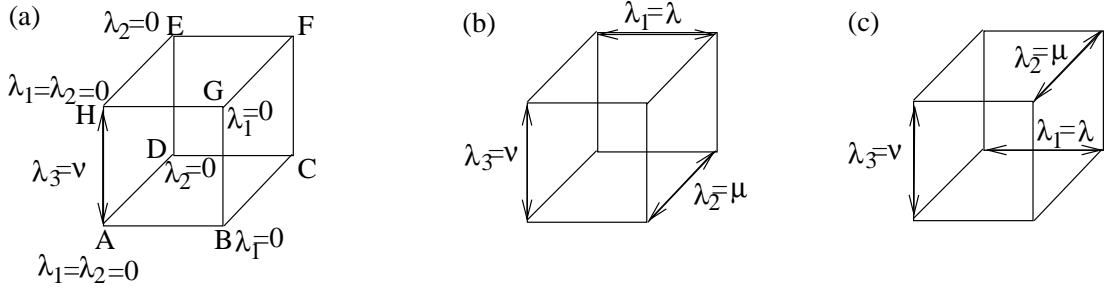}\]
 \caption{Flat connections of type (iii) in the cube: (a) Initial assumption,  (b) solution and (c) its mirror image as the only possible. In (b),(c) only the nonzero $\lambda_i$ are shown.}
 \end{figure}
 
Now suppose that $\lambda_2=\mu\ne 0$ at corners $B$ and $C$ (the two must be the same value) as shown by the arrowed edge in Figure~2(b). Then at $B$ we must have $\lambda_3=0$ to avoid a split (to avoid the existence of a point with two non-zero $\lambda_i$). In this case $\lambda_1=\lambda_3=0$ also at $C$. Hence $\lambda_1=0$ at $D$. We conclude also that $\lambda_3=0$ at $D$ (and hence all $\lambda_i=0$ at $D$) for if not,  we could deduce the same values at $E$ and hence that $\lambda_1=0$ at $F$, which would be a split with  $\lambda_1\equiv 0$.  Then $\lambda_1=\lambda\ne0$ at $E$ and $F$ (to avoid a split with $\lambda_1\equiv 0$). Hence $\lambda_2=\lambda_3=0$ at $E$ and $F$. Hence $\lambda_2=0$ also at $G$, and since $\lambda_3=0$ at $A$ it must also vanish at $G$, i.e. all three $\lambda_i$ vanish at $G$. The solution is then fully determined by the three non-zero values $\lambda,\mu,\nu$ and all three $\lambda_i$ vanishing at $D,G$ as shown in part (b) of the figure. We mark only the non-zero edges, which imply those value on their endpoints; all other values are zero.

Alternatively, if $\lambda_2=0$ at $B$ and $C$, then $\lambda_2=\mu\ne 0$ at $F$ and $G$ (the arrowed edge shown in Figure~2 part (c)) to avoid a split with $\lambda_2\equiv 0$. Hence $\lambda_3=0$ at $G$ to avoid a maximal solution, and hence also at $B$ (so all three $\lambda_i=0$ at $B$). Moreover, the values $\lambda_1=\lambda_3=0$ are transported to $F$. Hence $\lambda_3=0$ at $C$ also, and $\lambda_1=0$ at $E$ also. Finally, $\lambda_1=\lambda\ne0$ at $C,D$ (the final arrowed edge shown in part (c)) to avoid a split $\lambda_1\equiv0$. This transports $\lambda_3=0$ also to $D$ and hence to $E$, therefore we  deduce the mirror image solution to the one above, where the $\lambda_i=0$ now at $B,E$ as shown in part (c) of the figure.
  
The explicit formula in the first case, if $A$ is the origin of a standard cube,  is
 \[ \alpha=x_2x_3\lambda\tau_1+x_1(1-x_3)\mu\tau_2+(1-x_1)(1-x_2)\nu\tau_3-\vartheta,\quad\lambda,\mu,\nu\in \R_{>0}.\]
 Of course, we can rotate this solution by picking any other origin and initial non-zero edge, and we also  have the mirror image solution.  Finally, phases can be removed by gauge transformation in a similar manner to the above.  This  exhausts the moduli space for flat connections for the cube $n=3$ up to gauge equivalence. 
 
 \subsection*{Acknowledgements} The author would like to thank Sanjay Ramgoolam for posing the problem and for interesting discussions.


\begin{thebibliography}{0}\itemsep 0pt 

\bibitem{AkrMa:bra}
E. Akrami \& S. Majid.
\newblock Braided cyclic cohomology and nonassociative geometry.
\newblock {\em J. Math. Phys.} 45:3883--3911, 2005.

\bibitem{AlbMa:qua}
H. Albuquerque \& S. Majid.
\newblock Quasialgebra structure of the octonions.
\newblock {\em J. Algebra} 220:188-224, 1999.

\bibitem{AlbMa:Zn}
H. Albuquerque \& S. Majid.
\newblock $Z_n$-Quasialgebras.
\newblock {\em Textos de Mat. de Coimbra Ser. B} 19: 57-64, 1999.

\bibitem{AlbMa:cli}
H. Albuquerque \& S. Majid.
\newblock Clifford algebras obtained by twisting of group algebras.
\newblock {\em J. Pure Applied Algebra} 171:133-148, 2002.

\bibitem{BegMa:sem} E.J.  Beggs \& S. Majid.
\newblock Semiclassical differential structures.
\newblock Math.QA/0306273, to appear {\em Pac. J. Math.}

\bibitem{BegMa:qua} E.J. Beggs \& S. Majid.
\newblock Quantization by cochain twists and nonassociative differentials.
\newblock Preprint, math.QA/0506450.

\bibitem{BHTZ}
B. Bernevig, J. Hu, N. Toumbas \& S-C. Zhang.
\newblock Eight-dimensional quantum Hall effect and "octonions".
 \newblock {\em Phys. Rev. Lett.} 91, 236803, 2003.

 
\bibitem{Con}
A.~Connes.
\newblock {\em Noncommutative Geometry}.
\newblock Academic Press, 1994.

\bibitem{Dix}
G. Dixon.
\newblock  {\em Division Algebras: Octonions, Quaternions, Complex Numbers and the Algebraic Design of Physics}. 
\newblock Kluwer, Dordrecht, 1994.5   
  
  
\bibitem{Dri:qua}
V.G. Drinfeld.
\newblock Quasi{H}opf algebras.
\newblock {\em Leningrad Math. J.}, 1:1419--1457, 1990.

\bibitem{HoRam:geo}
P. M. Ho \& S. Ramgoolam.
\newblock Higher dimensional geometries from matrix brane constructions.
\newblock {\em Nucl. Phys. B} 627:266, 2002.

 
 \bibitem{Ma:dia}
 S. Majid.
 \newblock Diagrammatics of Braided Group Gauge Theory.
 \newblock {\em J. Knot Th. Ramif.} 8: 731-771, 1999.
 
 \bibitem{Ma:rem}
 S. Majid.
 \newblock Some remarks on quantum and braided group gauge theory.
 \newblock {\em Banach Center Publications} 40:335--349, 1997.
 
 \bibitem{Ma:book} S.~Majid.
\newblock {\em Foundations of Quantum Group Theory}.
\newblock Cambridge Univ. Press, 1995.

 

 \bibitem{MaRai:ele}
 S. Majid \& E. Raineri
 \newblock Electromagnetism and gauge theory on the permutation group $S_3$.
 \newblock {\em J. Geom. Phys.} 44:129-155, 2002.
 
 \bibitem{Ma:yan}S. Majid. 
 \newblock Noncommutative differentials and Yang-Mills on permutation groups $S_N$.
 \newblock  {\em Lect. Notes Pure Appl. Maths} 239:189-214, 2004.  Marcel Dekker.
 
 \bibitem{Ma:phy}
 S. Majid
 \newblock Noncommutative physics on Lie Algebras, $Z_2^n$ lattices and Clifford algebras.
 \newblock In {\em Clifford Algebras: Application to Mathematics, Physics, and Engineering}, ed. R. Ablamowicz, pp. 491-518.  Birkhauser, 2003.
 
 \bibitem{Ram:sph}
S. Ramgoolam.
\newblock On spherical harmonics for fuzzy spheres in diverse dimensions.
\newblock {\em  Nucl. Phys. B} 610:461, 2001.

 
 \bibitem{Ram:gau}
 S. Ramgoolam.
 \newblock Towards Gauge theory for a class of commutative and non-associative fuzzy spaces.
 \newblock Preprint hep-th/0310153.

 
\end{thebibliography}
\end{document}